\documentclass[12pt,  reqno]{amsart} 
\usepackage{latexsym} 
\usepackage{amssymb} 
\usepackage{amsmath}
\usepackage{dsfont}
\usepackage{enumerate}
\usepackage{graphicx}

\usepackage[utf8]{inputenc}

\usepackage[Conny]{fncychap}
\usepackage{eufrak}
\usepackage{hyperref}  

\usepackage{lipsum}

\usepackage{booktabs}
\usepackage{graphicx}
\newtheorem{theorem}{Theorem}[section]

\newtheorem{definition}{Definition} [section]
\newtheorem{lemma}{Lemma}[section]
\newtheorem{proposition}{Proposition}[section]
\newtheorem{corollary}{Corollary}[section]
\newtheorem{example}{Example}[section]

\newcommand{\be}{\begin{equation}} 
\newcommand{\ee}{\end{equation}} 
\newcommand{\bea}{\begin{eqnarray}} 
\newcommand{\eea}{\end{eqnarray}} 

\newcommand{\eeq}{\end{equation}} 

\newcommand{\eeqn}{\end{eqnarray}} 
\newcommand{\beaa}{\begin{eqnarray*}} 
\newcommand{\eeaa}{\end{eqnarray*}}

\def\complex{\mathop{\raise .45ex\hbox{${\bf\scriptstyle{|}}$} 
     \kern -0.40em {\rm \textstyle{C}}}\nolimits} 
\def\hilbert{\mathop{\raise .21ex\hbox{$\bigcirc$}}\kern -1.005em {\rm\textstyle{H}}} 

\renewcommand{\varrho}{{\ell}}

\newcommand{\nproof}{\noindent {\bf Proof:\ }}

\def\squarebox#1{\hbox to #1{\hfill\vbox to #1{\vfill}}} 
\newcommand{\nqed}{\hspace*{\fill} 
           \vbox{\hrule\hbox{\vrule\squarebox{.667em}\vrule}\hrule}\bigskip}

\title[Applications of the Tsirelson spectral measures to noise filtering]{Some applications of the Tsirelson spectral measures to noise filtering problems}

\author{R\'{e}mi Lassalle}

\address{\textsc{Ceremade}, Universit\'{e} Paris 9 (Dauphine), PSL, Place du Mar\'echal De Lattre De Tassigny, 75775 Paris Cedex 16, France}
\email{lassalle@ceremade.dauphine.fr}
\begin{document}

\begin{abstract}
This paper investigates applications of the Tsirelson spectral measures to noise filtering problems within the classical Wiener noise framework. 
We particularly focus on those among those measures associated to square integrable Radon-Nikodym derivatives of  laws of observation signals absolutely continuous with respect to the
Wiener measure. Explicit formulas are provided to compute the values of the Tsirelson spectral measure of specific sets of interests. Then, several applications to filtering 
problems are considered, notably to the so-called innovation problem of filtering.

\end{abstract}

\maketitle

\keywords{Stochastic analysis, Noise filtering, Tsirelson's spectral measure, Innovation noise.} \\
{\bf Mathematics Subject Classification :}
 60H30; 60H40; 60H99.


\section{Introduction}

This paper investigates practical applications of the so-called \textsc{Tsirelson} \textit{spectral measures} to address precise \textit{noise filtering} problems, on the somehow tangible ground of the classical  \textsc{Wiener} noise. Within the range of stochastic analysis, the corresponding \textsc{Tsirelson}'s \textit{spectrum of noise} has been recently qualified to be an  important notion within the lines  of \textsc{It\^{o}'}s theorems on chaos expansions and martingale representations theorems in \cite{Watanabe15}, which motivates the study unrolled in the present paper. Although the whole paper will be performed within a classical continuous framework, at this stage it seems to be worth to notice  that the discrete counterpart to those tools has recently received  further developments and applications related to other fields as \textit{percolation} (for instance, see the Lemma 5.3. of \cite{GARBAN}).  Recall that these mathematical objects have been  introduced by \textsc{Tsirelson} (p.274 of \cite{TsirelsonRev}, Definition 3d2 of \cite{TsirelsonSaintFlour}) within a systematic study of mathematical models of what may be interpreted as pure noises : up to a normalization constant, they  provide \textit{quantitative} instruments which describe accurately the statistical properties of some random compact subsets, which are associated to square integrable functionals over the underlying probability space of a noise. Their behavior is known to yield a distinction between \textit{classical noises} and \textit{non-classical noises} (see \cite{LEJANRAIMOND}, \cite{TsirelsonRev}, \cite{WARRENWAT}), depending on whether all the related random compact sets are almost surely finite, or not (see \cite{Watanabe15}, and the references therein). 
Thus, up to a normalization, the present paper will essentially focus on the quantitative description of some finite random compact sets associated to the classical \textsc{Wiener} noise. As it is well known, still in this case, given $f\in L^2(\mu_W)$, where $\mu_W$ denotes the so-called \textsc{Wiener} measure, the existence of the \textsc{Tsirelson} spectral measure $\rho_f$ of $\mu_W$ with respect to $f$ readily follows from  the so-called  \textsc{It\^{o}-Wiener} chaos expansion theorem  of $L^2(\mu_W)$, as it will be recalled. 

For the sake of clarity, in the whole paper, the time parameter $t$ will run through the whole closed unit interval $\mathbb{T}:=[0,1]$ (rather than $\mathbb{R}$), so that the noise will be conveniently described here by a standard Brownian motion $(B_t)_{t\in\mathbb{T}}$ defined on a complete probability space $(\Omega,\mathcal{A}, \mathbb{P})$ which is endowed with a filtration $(\mathcal{A}_t)_{t\in\mathbb{T}}$ which satisfies to the usual conditions (it is right continuous and complete, see \cite{PROTTER}). Subsequently, we handle  ${\bf B}$ as an (equivalence class of suitably measurable)  application(s) ${\bf B} : \omega \in \Omega \to C([0,1],\mathbb{R})$ which to a given $\omega\in \Omega$ associates the continuous curve $t\in \mathbb{T}\to B_t(\omega)\in  \mathbb{R}$. The \textit{noise filter} will be described by an (equivalence class of measurable) application(s) of the form $${\bf F}^u: {\bf x}\in C([0,1],\mathbb{R}) \to {\bf x} + {\bf u}({\bf x}) \in C([0,1],\mathbb{R}),$$ where $t\to {\bf u}_t({\bf x})$ is assumed to be absolutely continuous, with finite energy, and determined at each $t$ by the sole past values $\{{\bf x}(s) : s\in [0,t]\}$ of ${\bf x}\in C([0,1],\mathbb{R})$, outside some suitable negligible sets.  This filter turns the \textit{input} noise $(B_t)_{t\in\mathbb{T}}$ into an \textit{output observation} signal which we describe by a stochastic process $(Y_t)$, which in functional form reads $${\bf Y} = {\bf F}^u\circ {\bf B} ,$$ $\circ$ denoting the pullback of equivalence classes of maps (more accurately of morphisms of probability spaces, see \cite{AirMal}), which will require  to take into account precise aspects of absolute continuity of some deterministic \textit{causal transport plans} of $\mu_W$, and for short we shall also denote it by ${\bf F}^u({\bf B})$. Similarly, when ${\bf h}\in H^1$ is an element of the so-called \textsc{Cameron-Martin} space (see \cite{MALLIAVIN}), so that $t\in [0,1] \to {h}_t\in \mathbb{R}$ is an absolutely continuous function with square integrable derivative $(\dot{h})$, which describes a deterministic signal of finite energy, we denote by ${\bf Y}^h$ the output which is obtained by applying the filter to the noisy input signal ${\bf X}:= {\bf B} +{\bf h}$, which with a similar notation reads $${\bf Y}^h := {\bf F}^u({\bf B} +{\bf h}),$$ which is well defined. Under further assumptions, notably that the \textsc{Radon-Nikodym} derivative $f_Y:=\frac{dp_Y}{d\mu_W}$ of the absolutely continuous probability $p_Y$ with respect to $\mu_W$  is square integrable, this paper investigates applications of the \textsc{Tsirelson} spectral measure $\rho_{f_Y}$ to associated noise filtering problems. More accurately, this work is directed toward the following practical cases of study :

\begin{enumerate}[{\bf (1)}]
\item The question of the existence of a \textit{realizable inverse filter} described by some ${\bf H}^u$ to the real-time filter described by ${\bf F}^u$, and its application to the \textit{statistical estimation} of a deterministic input signal $(h_t)$ by a \textit{robust} estimator $( \hat{h}_t)$, for $t\in \mathbb{T}$. 

\item The so-called \textit{innovation problem} (see \cite{ALLINGERM}, \cite{MEYERF}). In a nutshell, this problem, which will be recalled accurately below, is to decide whether the data of a specific \textit{observable noise}, which is called the \textit{innovation} process, provides a sufficient summary to access in real-time to the whole information on past values of the observation process $(Y_t)$. 

\end{enumerate}

As far as $(1)$ is concerned, we obtain a necessary and sufficient condition which is based on \textsc{Tsirelson}'s spectral measures in Theorem~\ref{Theorem1TSMF}, for the existence of a realizable inverse filter. As an application, we establish a sufficient condition on $\rho_{f_{\bf Y}}$ for the existence of the real-time unbiased \textit{robust} statistical estimator $\widehat{h}_t:= { H}^{\bf u}_t({\bf Y^h})$ of the value at $t$ of the deterministic signal described by $h_t\in \mathbb{R}$. Here, we emphasize that the problematic of existence of such an inverse filter is a classical problem of the spectral analysis of time series. As a classical reference, we refer here to  \cite{KOOPMANS} p.105 where similar problems  are stated with linear filters and stationary random inputs : within the framework and hypothesis of \cite{KOOPMANS} the linearity and the stationarity allow to use classical spectral measures. However, within our  framework, the noise is $\textsc{Wiener}$ and the filter is not necessarily linear. Still, Theorem~\ref{Theorem1TSMF} shows that within the specific circumstances which we encounter within this context, the \textsc{Tsirelson} spectral measure $\rho_{f_Y}$ of $\mu_W$ with respect to $f_Y$ can be used as a somehow analogous substitute to the usual spectral measure to  provide a criterion for the existence for the inverse filter ${\bf H}^u$. The proof of the latter takes advantage of recent advances in stochastic analysis from works of \textsc{D.Feyel}, \textsc{A.S. \"{U}st\"{u}nel}, \textsc{M.Zakai} and their collaborators (for instance, see \cite{FEYEL1}, \cite{FEYEL2}, \cite{UZ07}, \cite{USTUNEL}) on distinct but related problems which, as a by product, can be used to study the invertibility of filters in particular cases. More accurately, within the proof a key role will be played by a \textit{causal} \textsc{Monge-Amp\`{e}re} equation, while some values of $\rho_{f_Y}$ plays a role somehow analogous to the relative entropy in \cite{USTUNEL}. Turning now to the case $(2)$, by following a parallel strategy, we show that the so-called \textit{innovation problem} of filtering can be also conveniently addressed by using the \textsc{Tsirelson} spectral measure of $f_Y$.

The structure of this paper is the following. In Section~\ref{Section2} we fix the notation for the whole paper. In particular, the definition of the \textsc{Tsirelson} spectral measures is recalled within the particular framework of the \textsc{Wiener} noise, under the specific hypothesis which we adopt here. Then, in Section~\ref{Section3}, under precise hypothesis,  we obtain explicit formulas for the values of the \textsc{Tsirelson} spectral measure $\rho_{f_Y}(A)$, for specific \textsc{Borel} subsets $A$ of the set $Comp(\mathbb{T})$ of compact subsets of $\mathbb{T}$, which is endowed with the corresponding \textsc{Hausdorff} distance. Those formulas are then applied in Section~\ref{Section4} to state Theorem~\ref{Theorem1TSMF} which corresponds to the case $(1)$ above. Then, Corollary~\ref{Corollary1} of Section~\ref{Section5} provides a precise result toward the application $(2)$ to the innovation problem. Finally, in Section~\ref{Section6} which also provides further perspectives, a probabilistic normalization is performed, so that the results of Section~\ref{Section4} first yield an explicit lower bond for the probability of a precise random set associated to $(Y_t)$ to be empty, and then provide a criterion for $(1$) which is founded on the law of this specific spectral random set.

\label{Section1}

\section{Preliminaries and notation}
\label{Section2} 

\subsection{Notation for stochastic processes}

Within the filtering model which we consider here, any realization of random signals on the interval of time $[0,1]$ will be described by a continuous function ${\bf \omega} : t \in [0,1] \to \omega(t)\in \mathbb{R}$. Therefore, we denote by $C([0,1], \mathbb{R})$  the space of continuous real valued functions defined on $[0,1]$, which we endow with the norm  $\|.\|_{\infty}$ of uniform convergence, so that it is turned into a separable \textsc{Banach} space, whose associated \textsc{Borel} sigma-field will be denoted by $\mathcal{B}(C([0,1],\mathbb{R}))$. Within this point of view, the law of a random continuous signal is described by an element $\nu$ of the set ${\bf M}_{1}(C([0,1], \mathbb{R}))$ of \textsc{Borel} probability measures on  $C([0,1], \mathbb{R})$. Among those laws, the standard \textsc{Wiener} measure $\mu_W \in{\bf M}_{1}(C([0,1], \mathbb{R}))$  will play a key role within this model, as it will be used to model a \textit{white noise}, since the latter is well known to be conveniently interpreted as some generalized derivative of the Brownian motion. Given a standard Brownian motion $(B_t)_{t\in[0,1]}$, with $B_0=0$, $\mathds{P}-a.s.$, which is defined on a complete probability space $(\Omega, \mathcal{A}, \mathds{P})$, and denoting by ${\bf B} : \Omega \to C([0,1],\mathbb{R})$ the application which to $\omega\in \Omega$ associates the curve $t\to B_t(\omega)$ of a continuous modification, outside some $\mathbb{P}-$ null set, recall that we have $\mu_W = p_B$, where $p_B:= {\bf B}_\star \mathds{P}$ denotes the pushforward of the probability measure $\mathds{P}$ by the $\mathcal{A}/ \mathcal{B}(C([0,1],\mathbb{R}))-$ measurable function ${\bf B}$, that is, $p_B(A):=\mathbb{P}(\{\omega \in \Omega : {\bf B(\omega)}\in A\})$, $\forall A\in \mathcal{B}(C([0,1],\mathbb{R}))$. Realizations of signals of finite energy will be described as elements of the classical \textsc{Cameron-Martin} space $H^1$ (see \cite{MALLIAVIN}), which is the linear subspace of $C([0,1], \mathbb{R})$ defined by $$H^1= \left\{ {\bf h}: [0,1] \to \mathbb{R}^d : \ {\bf h}=\int_0^. {\bf \dot{h}_s} ds \ , \ \int_0^1\|{\bf \dot{h}_s}\|^2_{\mathbb{R}^d}ds <+\infty \ \right\}; $$ recall that the latter takes its name from a celebrated sequence of papers of \textsc{R.H. Cameron}, and \textsc{W.T. Martin} (for instance see \cite{34.}, \cite{35.}).

  When it is endowed with the scalar product $<.,.>_{H^1} : (h,k)\in H^1\times H^1 \to \int_0^1 \dot{h}_s \dot{k}_s ds\in \mathbb{R}$, $H^1$ is a separable \textsc{Hilbert} space.  Subsequently, to handle random signal models as time evolves together with martingales, and to overcome technical aspects, we will require to work on  a complete stochastic basis  $(\Omega, \mathcal{A}, (\mathcal{A}_t)_{t\in[0,1]}, \mathds{P})$, where $(\Omega, \mathcal{A},\mathds{P})$ denotes  a $\mathbb{P}-$ complete probability space, and $(\mathcal{A}_t)_{t\in[0,1]}$ a filtration which satisfies to the usual conditions. Furthermore, when  $X :\Omega \to \mathbb{R}$ is an $\mathcal{A}/\mathcal{B}(\mathbb{R})-$ measurable function, where $\mathcal{B}(\mathbb{R})$ denotes the \textsc{Borel} sigma field on $\mathbb{R}$, we use the standard notation of the mathematical expectation $\mathbb{E}_{\mathds{P}}\left[X\right]:= \int_\Omega X(\omega) \mathds{P}(d\omega)$, if $E_{\mathbb{P}}[|X|]<+\infty$, in which case  $E_{\mathbb{P}}[X]\in \mathbb{R}$ and $X$ is said to be $\mathbb{P}$-integrable, or if $X\geq 0$, $\mathbb{P}-a.s.$, in which case $E_{\mathbb{P}}[X]\in [0,+\infty]$. At this stage it seems to be worth to mention that within the whole paper, we adopt the usual convention according to which the $\mathbb{P}-$ completion $\mathcal{G}^{\mathbb{P}}$ of a sub-sigma field $\mathcal{G}$ of $\mathcal{A}$ is defined to be the smallest sigma-field on $\Omega$ among those which contain both all the elements of $\mathcal{G}$ and all the elements of $\mathcal{N}^{\mathds{P}} :=\{ N\subset\Omega : \exists A\in \mathcal{A}, N\subset A, \mathbb{P}(A) =0\}$. To handle random processes of finite energy, we define the \textsc{Hilbert} space $L^2(\mathbb{P},H^1)$ to be the set which is constituted of ($\mathds{P}$- equivalence classes) of  $\mathcal{A}/ \mathcal{B}(C([0,1], \mathbb{R}))$- measurable functions $f : \Omega \to C([0,1],\mathbb{R})$ such that $E_{\mathbb{P}}[\|f\|_{H^1}^2] <+\infty$. Furthermore, in order to take into account the \textit{information flows} described by $(\mathcal{A}_t)_{t\in[0,1]}$ (see \cite{CONT}, Definition 2.11, p. 39), we will further consider the subset $L^2_a(\mathbb{P},H^1)$ of the ${\bf u}\in L^2(\mathbb{P}, H^1)$, such that the random variable $f_t : \omega \in \Omega \to f_t(\omega)\in \mathbb{R}$ is $\mathcal{A}_t-$ measurable, $\forall t\in [0,1]$, for any (and then all) $\mathcal{A}/ \mathcal{B}(C([0,1], \mathbb{R})- $measurable function $f : \Omega \to C([0,1],\mathbb{R})$ whose equivalence classes is ${\bf u}$, $f_t(\omega)$ denoting the value of the continuous function $f(\omega)\in C([0,1],\mathbb{R})$ at point $t\in [0,1]$, $\forall \omega\in \Omega$. Then, given $t\in[0,1]$, to shorten subsequent statements and proofs, it is convenient to define the operator $\pi_t : u\in L^2_a(\mathbb{P}, H^1) \to \int_0^. 1_{[0,t]}(\sigma) \dot{u}_\sigma d\sigma \in  L^2_a(\mathbb{P}, H^1)$.

 In the particular case where $(\Omega, \mathcal{A}, \mathds{P})$ coincides with the canonical probability space $( C([0,1],\mathbb{R}),$ $\mathcal{B}(C([0,1],$ $\mathbb{R})), \nu)$, for some $\nu \in {\bf M}_1(C([0,1],\mathbb{R}))$, in this paper, we systematically endow it with the usual augmentation $(\mathcal{F}_t^\nu)$ of the filtration generated by the evaluation process $(W_t)_{t\in[0,1]}$,  where $W_t :  \omega \in C([0,1],\mathbb{R})\to \omega(t) \in \mathbb{R}$ denotes the function which to any $\omega$ associates its value $W_t(\omega)= \omega(t)$ at $t\in[0,1]$ ; recall that $(\mathcal{F}_t^\nu)$ is right-continuous and $\nu-$ complete.  In this latter case, we simply denote by $L^2_a(\nu, H^1)$ for the corresponding $L^2_a(\mathds{P}, H^1)$ space, which doesn't seem to yield any confusion below. Aside filtrations, to perform computations, given a stochastic process $(Y_t)_{t\in[0,1]}$ defined on a complete probability space $(\Omega, \mathcal{A}, \mathds{P})$, and given $s,t\in[0,1]$ such that $s\leq t$, we shall denote by $(\mathcal{G}_{[s,t]}^Y)_{t\in[s,1]}$ the usual augmentation of the filtration $(\sigma(Y_v-Y_u : s\leq u\leq v\leq t))_{t\in[s,1]}$, and whenever $Y_t=W_t$, $\forall t\in [0,1]$, and the underlying probability space is $( C([0,1],\mathbb{R}), \mathcal{B}(C([0,1], \mathbb{R})), \nu)$, for some $\nu \in {\bf M}_1(C([0,1],\mathbb{R}))$, we use the specific notation $\mathcal{F}_{[s,t]}^\nu$ to denote $\mathcal{G}_{[s,t]}^Y$, so that $\mathcal{F}_t^\nu = \mathcal{F}_{[0,t]}^\nu$, $\forall t\in [0,1]$. In particular, with this notation, the independency of increments of the Brownian motion yields the continuous product structure, with in particular $$\mathcal{F}_1^{\mu_W}= \mathcal{F}_{[0,s]}^{\mu_W}\otimes \mathcal{F}_{[s,1]}^{\mu_W}, \forall s\in [0,1]$$ where the sigma-field in the r.h.t. denotes the smallest sigma-field which contains both all elements of $\mathcal{F}_{[0,s]}^{\mu_W}$ and elements of the $\mu_W-$ independent sigma-field $\mathcal{F}_{[s,1]}^{\mu_W}$.

\subsection{The Tsirelson spectral measures}

We denote by $Comp([0,1])$ the set which is constituted by the compact subsets of $[0,1]$, including the empty-set $\emptyset$, which is endowed with the so-called  \textsc{Hausdorff} \textit{distance} $d_{Haus} : Comp([0,1])\times Comp([0,1]) \to [0,+\infty]$, defined by $$ d_{Haus}(K,\widetilde{K}):= \inf\left(\left\{ \epsilon \in (0,+\infty) : K \subset \widetilde{K}^\epsilon \  \text{and} \  \widetilde{K} \subset K^\epsilon \right\}\right),$$ for any $K,\widetilde{K}\in Comp([0,1])$, where for any subset $D$ contained in $[0,1]$, and $\epsilon \in (0,+\infty)$, we use the notation $D^{\epsilon}:=\{ t \in [0,1] : \exists s \in D, \ |s-t|<\epsilon\}$.  We denote by $\mathcal{B}(Comp([0,1]))$ the associated \textsc{Borel} sigma-field on $Comp([0,1])$.    Within our specific framework, we recall the definition of the \textsc{Tsirelson} spectral measures,  in particular see section 3d4 of
 \cite{TsirelsonSaintFlour}. Although this definition usually holds for much general noise models with rather elaborated proofs, by following \cite{TsirelsonSaintFlour}, we emphasize that within our specific framework, its existence  follows from the so-called \textsc{It\^{o}-Wiener} chaotic representation of $L^2(\mu_W)$, as it will be recalled in this section.  For the sake of clarity denote by $\mathcal{E}:=\{ \cup_{i=1}^n [s_i, t_i] : 0\leq s_1\leq t_1 \leq...\leq s_s\leq t_n, \ n\in \mathbb{N}\}$, and for $E= \cup_{i=1}^n [s_i, t_i]\in \mathcal{E}$, $\mathcal{F}_E^{\mu_W}$ denotes the sigma-field generated by the $\mathcal{F}_{[s_i,t_i]}^{\mu_W}$, $\forall i\in \{1,...,n\}$.

\begin{definition} (Tsirelson's spectral measures within this specific framework.) \label{Tsimesdef} 
Given $f\in L^2(\mu_W)$ there exists a unique finite \textsc{Borel} measure $\rho_f$ on  $\text{Comp}([0,1])$ such that $$\rho_f(\{K\in \text{Comp}([0,1]) : K\subset E\}) = E_{\mu_W}\left[\left(E_{\mu_W}[f| \mathcal{F}_E\right])^2\right],$$ holds, $\forall E\in \mathcal{E}$. $\rho_{f}$ is called the \textsc{Tsirelson} spectral measure of $\mu_W$ with respect to $f$.
\end{definition}

As far as the existence is concerned, for the reader's convenience, we now recall how this measure can be explicitly defined for $f\in L^2(\mu_W)$ by setting \begin{equation} \rho_f(\mathcal{K}) :=  |f_0|^2 \mathds{1}_{\mathcal{K}}(\emptyset)+ \sum_{n=1}^{+\infty} n! \int_{[0,1]^n} |f_n(t_1,...,t_n)|^2 \mathds{1}_{\mathcal{K}}(\{t_1,...,t_n\}) \lambda_{\mathbb{R}^n}(dt_1,...,dt_n), \label{tsimot} \end{equation}  for any $\mathcal{K} \in \mathcal{B}(Comp([0,1]))$, where $\lambda_{\mathbb{R}^n}$ denotes the \textsc{Lebesgue} measure on $\mathbb{R}^n$, while $$f= \sum_{n=0}^{+\infty} I_n(f_n), \ \mu_W-a.s.$$  denotes the \textsc{It\^{o}-Wiener} chaos expansion (section XXI of \cite{DM5}, \cite{IKEDA}, \cite{MALLIAVIN}, p. 68 of \cite{MEYERQUANTUM}, \cite{NOURDIN}, \cite{PRIVAULT},  Theorem 3.7 p.202 of  \cite{YORBOOK}) of $f$ (the limit is in $L^2(\mu_W)$),  with $f_0=E_{\mu_W}[f]\in \mathbb{R}$, with $f_n\in L^2([0,1]^n)$ which is a symmetric  square integrable real-valued function defined on $[0,1]^n$ which is endowed with the restriction of the measure $\lambda_{\mathbb{R}^n}$ to the \textsc{Borel} $\sigma-$ field $\mathcal{B}([0,1]^n)$, for any $n\in \mathbb{N}$, and where  $$I_n(f_n):= \int_{[0,1]^n} f_n dW^{\otimes n}  := n! \int_0^1(\int_0^{s_1}...(\int_0^{s_{n-1}} f_n(s_1,...,s_n) dW_{s_n})...)dW_{s_1},$$ denotes iterated \textsc{It\^{o}}'s integrals, with $(n+1)! = (n+1)n!$ and $0!=1$,  while $I_0(f_0)(\omega) :=f_0$, $\mu_W-$ a.s. coincides with $f_0\in \mathbb{R}$ outside some $\mu_W-$ negligible set  ; notice that, since  the function $$(t_1,...,t_n)\in [0,1]^n \to \{t_1,...,t_n\} \in Comp([0,1])$$ is continuous and therefore \textsc{Borel} measurable, for any $n\in \mathbb{N}$, the r.h.t. of~(\ref{tsimot}) is well defined. Then, it  follows that $\rho_f$ meets the hypothesis of Definition~\ref{Tsimesdef}, while the uniqueness is clear.

\begin{example}\label{Example1}
 For $f= \frac{W_1^4}{\sqrt{E_{\mu_W}[W_1^8]}}$, $\rho_f$ is a probability measure. In this particular case, $\rho_f$ is the law of a random set, which is obtained from an i.i.d. sequence of random variables $(\xi_n)_{n\in \mathbb{N}}$ of uniform law on $]0,1[$ ($\xi_1\sim \mathcal{U}(]0,1[)$), which is defined on some probability space $(\Omega, \mathcal{A}, \mathbb{P})$, where is also defined an independent integer-valued random variable $N : \Omega \to \mathbb{R}$, which is distributed according to the law $p_N:= \sum_{n=0}^{+\infty} p_N(\{n\}) \delta_n$, where $\delta_n$ denotes the \textsc{Dirac} mass concentrated at $n\in \{0\}\cup \mathbb{N}$, and where $$p_N(\{n\}) := \rho_f(\{K\in Comp([0,1]) : Card(K) = n\}),$$ $\forall n\in \{0\}\cup \mathbb{N}.$ Then $\rho_f$ is the law of the random set $$C : \omega \in \Omega \to  \begin{cases}   \{\xi_1(\omega),..., \xi_{N(\omega)}(\omega)\} \in Comp([0,1])
        & \text{if } \  N(\omega) \geq 1 \\
       \emptyset \in Comp([0,1]) &  \text{if} \  N(\omega) = 0
       \end{cases} .$$ 
\end{example}

\begin{table}[htbp] \label{Table11}
  \centering
 \caption{Numerical table : numerical approximations with $10^{-5}$ accuracy of values of the function $n\to p_N(\{n\})$ of Example~\ref{Example1}.}

    \begin{tabular}{lcccccc}
        \hline
        \hline
          $n$                   & 0        & 1          & 2        & 3  & 4   &  $>$4       \\
          $p_N(\{n\})$               & 0.08571          & 0          &  0.68571          & 0  & 0.22857    & 0      \\
        \hline
    \end{tabular}
     
  \label{tab:grav}
\end{table}

\begin{example}\label{Example2}
For $\beta\in(0,+\infty)$, set $f:= \exp(\sqrt{\beta} W_1 -\beta)$, so that once again $\rho_f$ is a probability since $E_{\mu_W}[f^2]=1$. Then, the probability law $p_N$ defined by $p_N(\{n\}) := \rho_f(\{K\in Comp([0,1]) : Card(K) = n\})$, $\forall n\in \mathbb{N}$, is now a \textsc{Poisson} distribution of parameter $\beta$. 
\end{example}

\begin{figure}
   \center 
    \includegraphics[width=4cm]{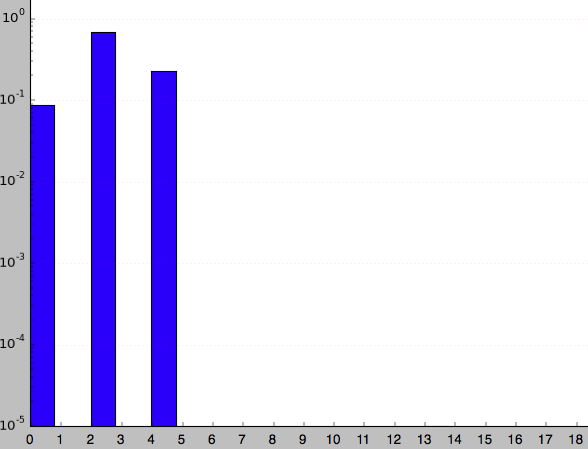}
      \includegraphics[width=4cm]{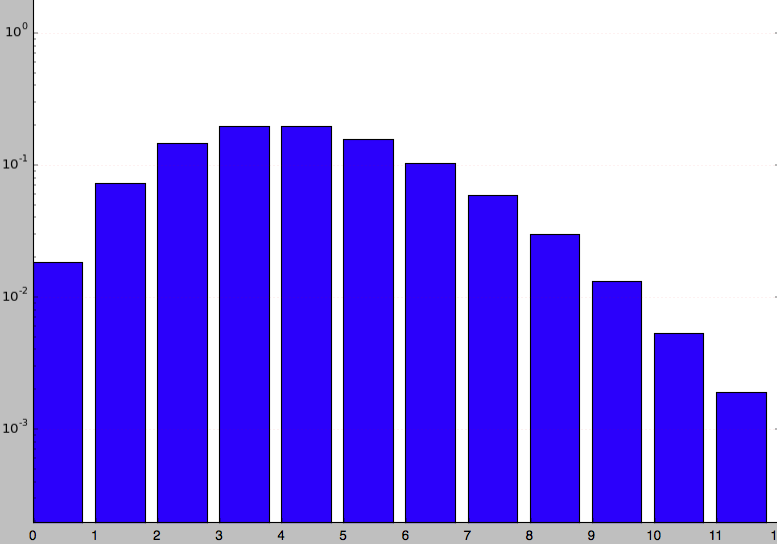}
    \caption{The figure on the left (resp. on the right) illustrates Example~\ref{Example1} (resp. Example~\ref{Example2}). It represents a restriction of the graph of  the \textsc{Tsirelson} spectral mass function $n  \in \{0\}\cup \mathbb{N} \to p_N(\{n\}) \in [0,1]$, with the data of Table~1 (resp. by taking $\beta= 4$ in Example~\ref{Example2}).}
   \end{figure}

\subsection{The noise filter}

In the whole paper, $(B_t)_{t\in[0,1]}$ will denote a standard brownian motion defined on a complete stochastic basis  $(\Omega, \mathcal{A}, (\mathcal{A}_t)_{t\in[0,1]}, \mathds{P})$, 
where $(\mathcal{A}_t)$ satisfies to the usual conditions. Given ${\bf h}\in H^1$, we may consider noisy input signals of the form $$X_t = { B}_t + h_t, \ \forall t\in [0,1], \mathbb{P}-
a.s..$$ It is a well known consequence of the so-called \textsc{Cameron-Martin-Segal} theorem (see \cite{MALLIAVIN}) that the law $p_X$ is actually equivalent to the \textsc{Wiener} 
measure $p_B=\mu_W$, which we denote by $p_X\sim \mu_W$.  We further assume that this noisy input signal is then turned into the observation signal $(Y_t)$ by a not necessarily 
linear filter ${\bf F}^{\bf u}$  which is determined by a given  ${\bf u}\in L^2_a(\mu_W,H^1)$,  from the formula  \begin{equation} {\bf F}^{\bf u} := \mathds{I}_{C([0,1],\mathbb{R})} +{\bf 
u}  \label{FUDEF}, \end{equation} so that we define \begin{equation}{\bf Y^h}={\bf F}^u\circ {\bf X},\label{YDEF00}\end{equation} $\circ$ denoting the usual pullback of 
(equivalence classes) of measurable functions which is well defined here, since $p_X\sim \mu_W$ (i.e. equivalent), and for short we may also denote it by $ {\bf F}^u({\bf X})$. The 
latter means that $$Y_t^{\bf h} = {F}^{\bf u}_t\circ {\bf X}, \ \forall t\in [0,1], \mathds{P}-a.s.,$$ where ${ F}^{\bf u}_t = W_t\circ {\bf F}^{\bf u}$,  and subsequently, to simplify the 
notation, we may also denote it by ${Y_t^{\bf h}}={ F}^{\bf u}_t({\bf X})$, $\forall t\in [0,1]$, while  in the particular case where ${\bf h}={\bf 0_{H^1}}$ we shall denote ${\bf Y^{0_{H^1}}}
$ by ${\bf Y}$. Notice here that,  for given $t\in [0,1]$, the value of ${ F}^{\bf u}_t({\bf X}):={ F}^{\bf u}_t\circ {\bf X}$ may be determined by the whole past values $\{ X_s, s \in[0,t]\}$ of 
the input $(X_s)_{s\in[0,1]}$, in a not necessarily linear fashion. Although filters of the form~(\ref{YDEF00}) may seem to be quite restrictive with respect to some of the forthcoming 
computations, it turns out to be taylor made to obtain a  compact statement of our main results Theorem~\ref{Theorem1TSMF}, and Corollary~\ref{Corollary1} on the innovation 
conjecture, while it emphasizes the clarity of the method. Sunsequently, $p_{\bf Y}$ (resp. $p_{{\bf X}}$) denotes the law of ${\bf Y}$ (resp. ${\bf X}$) that is, the \textsc{Borel} 
probability measure on $C([0,1], \mathbb{R})$ which is the pushforward of   $\mathbb{P}$ by the $C([0,1], \mathbb{R})-$ valued (equivalence class of) measurable function(s) ${\bf Y}
$ (resp. ${\bf X}$) which is defined on $(\Omega, \mathcal{A},\mathbb{P})$. In the particular case where ${\bf X} ={\bf B}$, notice again that we have $p_{\bf B} = \mu_W$.

\subsection{The innovation conjecture}

In this whole subsection we focus on the case where $h_t=0$, $\forall t\in [0,1]$, so that the input signal $(X_t)$ coincides with the noise $(B_t)$ ; as it will be clear below, within our 
specific purposes, the case where ${\bf h}\neq {\bf 0_{H^1}}$ can actually be reduced to this case, by considering the composition of the filter described by ${\bf F}^u$ with the filter 
described by the function ${\bf T}_h := \mathds{I}_{C([0,1],\mathbb{R})} +{\bf h}$.  Hence, in this section, given ${\bf u}\in L^2_a(\mathbb{P}, H^1)$ we still define ${\bf Y} := {\bf 
Y^{0_{H^1}}}$ by~(\ref{YDEF}) with ${\bf F}^u$ given by~(\ref{FUDEF}) with now ${\bf X}= {\bf B}$, so that we have \begin{equation} {\bf Y } = {\bf F}^u\circ {\bf B} = {\bf B} + {\bf 
\widetilde{u}}, \label{YDEF} \end{equation} where ${\bf \widetilde{u}}\in L^2_a(\mathbb{P},H^1)$ is defined by $${\bf \widetilde{u}} := {\bf u}\circ {\bf B} = {\bf Y} - {\bf B}.$$

Recall that, in this case, the innovation process of $(Y_t)$ which we denote here by $(B_t^{\bf \widetilde{u}})$ is the $(\mathcal{G}_{[0,t]}^Y)-$ brownian motion on $(\Omega, \mathcal{A}, \mathds{P})$  (see \cite{BAIN}, \cite{26.}, \cite{MEYERF}) which is defined by $$B_t^{\bf \widetilde{u}} := B_t + \int_0^t(\dot{\widetilde{u}}_s- E_{\mathds{P}}\left[\dot{\widetilde{u}}_s| \mathcal{G}_{[0,s]}^Y\right])ds, \ \forall t\in [0,1], \mathds{P}-a.s.,$$ so that it satisfies to $$Y_t = B_t^{\bf \widetilde{u}} + \int_0^t E_{\mathds{P}}\left[\dot{\widetilde{u}}_s| \mathcal{G}_{[0,s]}^Y\right]ds, \ \forall t\in[0,1], \mathds{P}-a.s..$$ At this stage, it seems to be worth to notice that the small increments of the so-called \textit{innovation noise process} $(B^{\bf \widetilde{u}}_t)$ may be interpreted as the part of increments of $(Y_t)$ which are not determined by the past of $(Y_t)$ and are independent to it under $\mathbb{P}$. A natural question is then to determine whether $(Y_t)$ is totally determined by its innovation $(B_t^{\bf \widetilde{u}})$, which is the so-called \textit{innovation problem} of \textsc{Frost} and \textsc{Kailath} (see \cite{ALLINGERM}, \cite{MEYERF} and the references therein). Subsequently, given $t\in[0,1]$,   $(Y_s)$ is said to satisfy to the \textit{innovation conjecture} of filtering on $[0,t]$, if and only if, the answer to the previous question is affirmative, that is, if and only if, the following coincidence of filtrations $$(\mathcal{G}_{[0,s]}^{{\bf B^{\bf \widetilde{u}}}})_{s\in[0,t]} = (\mathcal{G}_{[0,s]}^{\bf Y})_{s\in[0,t]}\  (\text{innovation conjecture  on  [0,t]}),$$ holds.  In cases, this problems may involves sharp aspects of filtrations theory (for instance, see \cite{Tsirelinvov}).

Denoting by $p_Y:={\bf Y}_\star \mathds{P}\in {\bf M}_1(C([0,1],\mathbb{R}))$ the law of $({Y}_t)$, it is known (see \cite{KZ}) to be absolutely continuous to the \textsc{Wiener} measure $\mu_W=p_B$ (\cite{WIENER}, see also \cite{KUO}),  so that the corresponding \textsc{Radon-Nikodym} derivative $\frac{dp_Y}{dp_B}\in L^1(p_B)$ is well defined. Under the further assumption that  $\frac{dp_Y}{dp_B}\in L^2(p_B)$, we can consider the \textsc{Tsirelson} spectral measure $\rho_{f_Y}$ with $f_Y:=\frac{dp_Y}{dp_B}$.  Subsequently, it is shown that under suitable hypothesis, such $(Y_t)$ which satisfy to the innovation conjecture can be characterized through $\rho_{f_Y}$. Finally, given ${\bf {\widetilde u}} \in L^2_a(\mathds{P}, H^1)$, $\delta^B{\bf \widetilde{u}}:=\int_0^1 \dot{\widetilde{u}}_s dB_s$ denotes the \textsc{It\^{o}} stochastic integral (\cite{44.}), while \begin{equation} \mathcal{E}\left(- \delta^B{\bf \widetilde{u}} \right):=\exp\left(-\delta^B{\bf \widetilde{u}} -\frac{\|{\bf \widetilde{u}}\|_{H^1}^2}{2} \right) \label{DoleansDadedef} \end{equation} denotes the so-called \textsc{Dol\'{e}ans-Dade} exponential (see \cite{PROTTER}).

\section{\textsc{Tsirelson}'s spectral measures associated to the \textsc{Radon-Nikodym} derivatives of $p_Y$ with respect to the law of the \textsc{Wiener} noise}
\label{Section3}

\begin{proposition} \label{Prop1}
Let $s,t\in [0,1]$ be such that $s\leq t$, and let  ${\bf u} \in L^2_a(\mu_W, H^1)$ which satisfies to the hypothesis \begin{equation} E_{\mu_W}\left[\exp\left(-\int_0^t \dot{u}_s dW_s -\frac{1}{2} \int_0^t \dot{u}_s^2 ds\right)\right] =1 \label{INTNOMPE2} .\end{equation} Further define ${\bf Y }$ by~(\ref{YDEF}), with ${\bf F}^{\bf u}$ given by~(\ref{FUDEF}), and set ${\bf \widetilde{u}}:={\bf Y}-{\bf B}\in L^2_a(\mathds{P}, H^1)$ so that ${\bf Y}= {\bf B} +{\bf \widetilde{u}}$. Then, we have
$$E_{\mathds{P}}\left[\mathcal{E}\left(- \delta^B \pi_t{\bf \widetilde{u}}  \right)  \middle|   \mathcal{G}_{[s,t]}^{\bf Y} \right]       =    \exp\left( -\int_s^t E_{\mathbb{P}}\left[\dot{\widetilde{u}}_\sigma| \mathcal{G}_{[s,\sigma]}^{\bf Y}  \right] dY_\sigma  + \frac{1}{2} \int_s^t \left(E_{\mathbb{P}}\left[\dot{\widetilde{u}}_\sigma| \mathcal{G}_{[s,\sigma]}^{\bf Y}  \right] \right)^2 d\sigma  \right),$$ where the notation of the l.h.t is given by~(\ref{DoleansDadedef}), while the integral in the r.h.t. denotes the \textsc{It\^{o}} stochastic integral with respect to the semi-martingale $(Y_{\widetilde{s}})_{\widetilde{s}\in[0,1]}$.
\end{proposition}
\nproof
Let $s,t\in[0,1]$ be such that $s\leq t$. For $n\in \mathbb{N}$, first define the $(\mathcal{G}_{[s,\sigma]}^{\bf Y})_{\sigma\in [s,t]}-$ stopping time  $$\tau_n:= \inf\left(\left\{\widetilde{s} \in 
[s,t] :  \int_s^{\widetilde{s}} \left(E_{\mathbb{P}}\left[\dot{\widetilde{u}}_\sigma \middle| \mathcal{G}_{[s,\sigma]}^{\bf Y}  \right] \right)^2 d\sigma>n \right\}\right)\wedge t,$$ with the 
convention $\inf(\emptyset) = +\infty$, and further notice that since $\tau_n$ is $\mathcal{G}_{[s,t]}^{\bf Y}-$ measurable, it can be written of the form $\tau_n = \sigma_n \circ {\bf Y}
$, $\mathbb{P}-a.s.$, where $\sigma_n$ is a $\left(\mathcal{F}_{[s,\widetilde{s}]}^{p_{\bf Y}}\right)_{\widetilde{s}\in[s,t]}-$ stopping time on the probability space ($C([0,1], \mathbb{R})$, $\mathcal{B}(C([0,1],
\mathbb{R}))^{p_{\bf Y}}$, $p_{\bf Y}$). Since ${\bf \widetilde{u}}\in L^2(\mathbb{P}, H^1)$,  the conditional \textsc{Jensen} inequality yields $\lim_{n\to +\infty} \tau_n =t$, $\mathbb{P}-
a.s.$.  Let $\theta_{s,t}$ be a continuous and bounded $\mathcal{F}_{[s,t]}^{p_{\bf Y}}-$ measurable real-valued function which is defined on $C([0,1],\mathbb{R})$, and define $
\theta^n_{s,t}:= \theta_{s,t}\circ \rho_n$ where $\rho_n : \omega \in  C([0,1], \mathbb{R}) \to \omega_{.\wedge \sigma_n(\omega)}\in C([0,1], \mathbb{R}) $. First notice that due 
to~(\ref{INTNOMPE2}) the 
\textsc{Cameron-Martin-Girsanov} theorem (\cite{GIRSANOV}) ensures that $(Y_{\widetilde{s}}- Y_s)_{\widetilde{s} \in [s,t]}$ is an $(\mathcal{A}_{\widetilde{s}})_{\widetilde{s} \in [s,t]}-$ Brownian motion starting from $0$ at  $s$ on the probability space $(\Omega, \mathcal{A}, \mathds{P}_u)$, where $\mathds{P}_u$ denotes the probability absolutely continuous (and also due to~(\ref{INTNOMPE2}) equivalent) with respect to $\mathbb{P}$ whose 
\textsc{Radon-Nikodym} derivative is defined by $$\frac{d\mathds{P}_u}{d\mathds{P}}= \mathcal{E}\left(- \delta^{\bf B} \pi_t {\bf \widetilde{u}}  \right), \  \mathbb{P}-a.s..$$ Therefore, we first obtain \begin{equation}E_{\mathds{P}}\left[\theta_{s,t}^n\circ {\bf Y} \mathcal{E}\left(- \delta^{\bf B} \pi_t {\bf  \widetilde{u}}  \right)\right] = E_{\mathds{P}}\left[\theta_{s,t}^n\circ {\bf B}\right].  \label{eqmplsl1}\end{equation} On the other hand, notice 
that \begin{equation} \label{ysdeflem} Y_{\widetilde{s}}-Y_s = N_{\widetilde{s}} +\int_{s}^{\widetilde{s}}E_{\mathds{P}}\left[\dot{\widetilde{u}}_\sigma| 
\mathcal{G}_{[s,\sigma]}^{\bf Y} \right] d\sigma, \ \forall \widetilde{s}\in [s,t], \mathbb{P}-a.s.,\end{equation} where we have defined 
\begin{equation} N_{\widetilde{s}} := B_{\widetilde{s}}-B_s + \int_{s}^{\widetilde{s}} \left(\dot{\widetilde{u}}_\sigma -E_{\mathds{P}}\left[\dot{\widetilde{u}}_\sigma| 
\mathcal{G}_{[s,\sigma]}^{\bf Y} \right]\right) d\sigma, \ \forall \widetilde{s}\in[s,t], \ \mathbb{P}-a.s.. \label{Ndefmpr} \end{equation} 
While~(\ref{ysdeflem}) ensures that $(N_{\widetilde{s}})_{\widetilde{s} \in[s,t]}$ is $(\mathcal{G}_{[s,\widetilde{s}]}^{\bf Y})_{\widetilde{s}\in[s,
t]}-$ adapted with $N_{s}=0$, $\mathbb{P}-a.s.$, since for $\widetilde{s}\in [s,t]$, we have $\mathcal{G}_{[s,\widetilde{s}]}^{\bf Y}\subset 
\mathcal{A}_{\widetilde{s}}$, from~(\ref{Ndefmpr}), the tower property yields $$ E_{\mathbb{P}}\left[N_{\widetilde{t}}-N_{\widetilde{s}}| 
\mathcal{G}_{[s,\widetilde{s}]}^{\bf Y} \right] =0, \ \mathbb{P}-a.s.,$$ so that from \textsc{Paul L\'{e}vy}'s criteria we conclude that $(N_{\widetilde{s}})_{\widetilde{s} \in[s,
t]}$ is a $(\mathcal{G}_{[s,\widetilde{s}]}^Y)_{\widetilde{s}\in[s,t]}$ Brownian motion starting from $0$ at time $s$ on $(\Omega, \mathcal{A},\mathbb{P})$. Henceforth, and until the end of the proof, for $n\in\mathbb{N}$, let $(M^n_{\widetilde{s}})_{\widetilde{s}\in[s,t]}$ be the stochastic process
defined on $(\Omega, \mathcal{A}, \mathbb{P})$ by $$M_{\widetilde{s}}^n:= \exp\left( -\int_s^{\widetilde{s}} 1_{[0,\tau_n]}(\sigma) E_{\mathbb{P}}\left[\dot{\widetilde{u}}_\sigma| \mathcal{G}_{[s,\sigma]}^{\bf Y}  \right] dN_
\sigma  - \frac{1}{2} \int_s^{\widetilde{s}} 1_{[0,\tau_n]}(\sigma) \left(E_{\mathbb{P}}\left[\dot{\widetilde{u}}_\sigma \middle| \mathcal{G}_{[s,\sigma]}^{\bf Y}  \right] \right)^2 d\sigma  \right),$$ $\forall {\widetilde{s}}\in [s,t], \mathbb{P}-a.s.,$ which is a non-negative $(\mathcal{G}_{[s,\widetilde{s}]}^{\bf Y})_{\widetilde{s}\in[s,t]}-$ martingale due to the \textsc{Novikov} criterion. For $\theta_{s,t}^n$ as above, from the \textsc{Cameron-Martin-Girsanov} theorem,  now applied to the $(\mathcal{G}_{[s,\widetilde{s}]})_{\widetilde{s} \in [s,t]}-$ Brownian motion $(N_{\widetilde{s}})_{\widetilde{s}\in [s,t]}$, it therefore follows that $$E_{\mathbb{P}}[M_{t}^n \theta_{s,t}^n\circ {\bf Y}] =E_{\mathbb{P}}[\theta_{s,t}^n\circ {\bf B}],$$
which together with~(\ref{eqmplsl1}) yields 
\begin{equation}E_{\mathbb{P}}[M_{t}^n \theta_{s,t}^n\circ {\bf Y}] = E_{\mathds{P}}\left[\theta_{s,t}^n\circ {\bf Y} \mathcal{E}\left(- 
\delta^{\bf B}{\bf\widetilde{u}}\right)\right].\label{MREDUCEDEF}\end{equation} Now, notice that the optional sampling theorem (see \cite{IKEDA}, Theorem 6.11, p.34) ensures that $(M_{t}^n)_{n\in \mathbb{N}}$ is a uniformly integrable $(\mathcal{G}_{[s, 
\tau_n]}^Y)_{n\in\mathbb{N}}-$ martingale with respect to the probability $\mathbb{P}$, $\mathcal{G}_{[s, 
\tau_n]}^Y$ denoting the corresponding $\sigma-$ field of the $(\mathcal{G}_{[s,\widetilde{s}]}^{\bf Y})_{\widetilde{s}\in [s,t]}-$ stopping time $\tau_n$ (see Definition 5.3. of \cite{IKEDA}, p.22), $\forall n\in \mathbb{N}$,  and that it is closed by the $\mathbb{P}-$ integrable random variable 
$$M_t:= \exp\left( -\int_s^{t}  E_{\mathbb{P}}\left[\dot{\widetilde{u}}_\sigma| \mathcal{G}_{[s,\sigma]}^{\bf Y}  \right] dN_\sigma  - \frac{1}{2} \int_s^{t} \left(E_{\mathbb{P}}\left[\dot{\widetilde{u}}_\sigma \middle| \mathcal{G}_{[s,\sigma]}^{\bf Y}  \right] \right)^2 d\sigma  \right).$$ Hence, by noticing that $t\in [0,1]$ is fixed,  the martingale convergence theorem (see Theorem 10 and Theorem 13 p.8-9 of \cite{PROTTER}) ensures that the discrete time martingale
$(M_{t}^n)_{n\in \mathbb{N}}$ converges $\mathbb{P}-$ almost surely to $M_t\in L^1(\mathbb{P})$, while we have $\rho_n({\bf Y})= Y_{.\wedge \tau_n}$, $\forall n\in \mathbb{N}$, where $\rho_n({\bf Y}):= \rho_n \circ {\bf Y}$, $\forall n\in \mathbb{N}$. The latter entails that $(\rho_n({\bf Y}))_{n\in \mathbb{N}}$ converges $\mathbb{P}-a.s.$ to ${\bf Y}$, and then, by continuity, that $(\theta_{s,t}^n({\bf Y}))_{n\in \mathbb{N}}$ converges almost surely to $\theta_{s,t}
({\bf Y})$, where $\theta_{s,t}^n({\bf Y}):=  \theta_{s,t}^n \circ {\bf Y}$, $\forall n\in \mathbb{N}$  (respectively $\theta_{s,t}({\bf Y}):=  \theta_{s,t}\circ {\bf Y}$). Since $M_t\in L^1(\mathbb{P})$, and since $(\theta_{s,t}^n)$ is uniformly bounded, the \textsc{Lebesgue} convergence theorem (see \cite{DM1}) yields

 \begin{equation} \lim_{n\to +\infty}E_{\mathbb{P}}[M_{t}^n \theta_{s,t}^n\circ {\bf Y}] = E_{\mathbb{P}}\left[M_t \theta_{s,t}\circ {\bf Y}\right] \label{rhtmpls},  \end{equation}
while \begin{equation} 
 \lim_{n\to +\infty} E_{\mathds{P}}\left[\theta_{s,t}^n\circ {\bf Y} \mathcal{E}\left(- \delta^{\bf B} {\bf  \pi_t \widetilde{u}}  \right)\right] = E_{\mathbb{P}}\left[\theta_{s,t}\circ {\bf Y} 
\mathcal{E}\left(- \delta^{\bf B} {\bf  \widetilde{u}}  \right)\right] \label{rhtmpls2} \end{equation}  follows similarly.

 Then, taking the limit in~(\ref{MREDUCEDEF}), from~(\ref{rhtmpls}) 
and~(\ref{rhtmpls2}), it follows that  $$ E_{\mathbb{P}}\left[\theta_{s,t}\circ {\bf Y} \mathcal{E}\left(- \delta^{\bf B} \pi_t{\bf  \widetilde{u}}  \right)\right] =E_{\mathbb{P}}\left[M_t \theta_{s,t}\circ{\bf Y}
\right],$$ so that, taking into account that $M_t$ is $\mathcal{G}_{[s,t]}^{\bf Y}-$ measurable, we get the result. \nqed

\begin{figure}
  \center 
    \includegraphics[width=6cm]{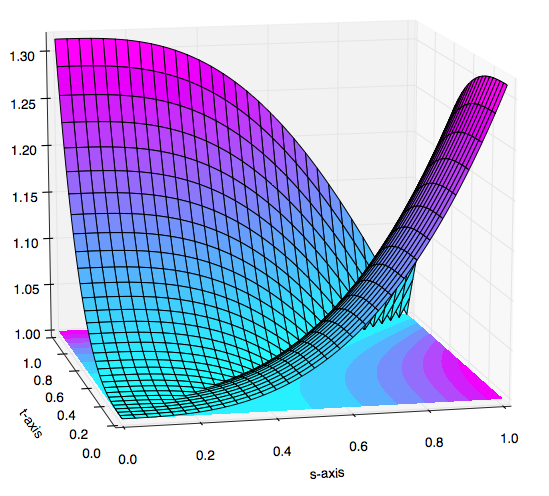}
    \caption{This figure illustrates Lemma~\ref{Lemme1}. It represents the graph of the function 
     $$(s,t)\in [0,1]\times [0,1] \to \rho_{f_Y}(\{K\in \text{Comp}([0,1]) : K \subset [s\wedge t, s \vee t]\}) \in \mathbb{R}_+,$$ for the particular case where $u=\int_0^. \sin(s)ds$, $\mu_W-a.s.$ ; $s\wedge t$ (resp. $s \vee t$) denotes $\min(s,t)$ (resp. $\max(s,t)$), $\forall s,t \in [0,1]$.}
 \end{figure}

\begin{lemma} \label{Lemme1} Given ${\bf u}\in L^2_a(\mu_W,H)$, further define ${\bf Y }$ by~(\ref{YDEF}), where ${\bf F}^{\bf u}$ is given by~(\ref{FUDEF}), and set ${\bf \widetilde{u}}:={\bf Y}-{\bf B}\in L^2_a(\mathds{P}, H^1)$. Further assume that $f_{\bf Y} \in L^2(\mu_W)$, where $f_{\bf Y}= \frac{dp_{\bf Y}}{dp_{\bf B}}$ denotes the \textsc{Radon-Nikodym} derivative of the absolutely continuous probability $p_{\bf Y}$ with respect to $p_{\bf B}=\mu_W$. Then, for any $0\leq s\leq t\leq 1$ such that~(\ref{INTNOMPE2}) holds with $t$,
 we have $$\rho_{f_Y}(J([s,t])) = E_{\mathbb{P}}\left[ \exp\left( \int_s^t E_{\mathbb{P}}\left[\dot{\widetilde{u}}_\sigma| \mathcal{G}_{[s,\sigma]}^Y  \right] dY_\sigma  - \frac{1}{2} \int_s^t \left(E_{\mathbb{P}}\left[\dot{\widetilde{u}}_\sigma| \mathcal{G}_{[s,\sigma]}^{\bf Y}  \right] \right)^2 d\sigma  \right)\right],$$ where \begin{equation} \label{JDEFPME} J([s,t]):= \{K \in \text{Comp}([0,1]) : K\subset [s,t]\} \in \mathcal{B}(\text{Comp}([0,1])),\end{equation} and where $\rho_{f_{\bf Y}}$ denotes the \textsc{Tsirelson} spectral measure (see Definition~\ref{Tsimesdef}) of $\mu_W$ with respect to $f_{\bf Y}$.
\end{lemma}
\nproof
From the definitions, we obtain \begin{equation} \rho_{f_{\bf Y}}(J([s,t])) = E_{\mathbb{P}}\left[ E_{{p_{\bf B}}}\left[\frac{dp_{\bf Y}}{d p_{\bf B}} \middle| \mathcal{F}_{[s,t]}^{{p_{\bf B}}}\right]\circ {\bf Y} \right]. \label{rhoydmts} \end{equation}  Since for any real valued functions $\theta_{s,t}$ which is bounded and $\mathcal{F}_{[s,t]}^{p_{\bf Y}}-$ measurable, the \textsc{Cameron-Martin-Girsanov} theorem, whose application is ensured by~(\ref{INTNOMPE2}),  yields $$E_{\mathbb{P}}\left[\theta_{s,t}\circ {\bf Y} E_{{p_{\bf B}}}\left[\frac{dp_{\bf Y}}{dp_{\bf B}}\middle| \mathcal{F}_{s,t}^{{p_{\bf B}}}\right]\circ {\bf Y} \mathcal{E}(-\delta^{\bf B} \pi_t{\bf \widetilde{u}})\right] = E_{\mathbb{P}}\left[\theta_{s,t}\circ {\bf Y}\right],$$ we obtain \begin{equation}  E_{{p_{\bf B}}}\left[\frac{dp_{\bf Y}}{dp_{\bf B}}\middle| \mathcal{F}_{[s,t]}^{{p_{\bf B}}}\right]\circ {\bf Y} E_{\mathbb{P}}\left[\mathcal{E}(-\delta^{\bf B} \pi_t {\bf \widetilde{u}}) \middle| \mathcal{G}_{[s,t]}^{\bf Y}\right] = 1, \ \mathbb{P}-a.s. . \label{cmapong} \end{equation} Substituting~(\ref{cmapong}) into (\ref{rhoydmts}), the result follows from Proposition~\ref{Prop1}. \nqed

\begin{remark} \label{remark1M}
By taking $s=0$,~(\ref{cmapong}) boils down to the \textit{causal} \textsc{Monge-Amp\`{e}re} equation as Proposition~3.4. of  \cite{FEYEL2}.
\end{remark}

\section{First application : realizable inverse of a real-time noise filter}
\label{Section4}

\begin{theorem} \label{Theorem1TSMF}
Let $t\in [0,1]$ and let ${\bf u}\in L^2_a(\mu_W, H^1)$ be such that~(\ref{INTNOMPE2}) holds.  Further define ${\bf Y}$ by~(\ref{YDEF}),  set ${\bf \widetilde{u}}:={\bf Y} -{\bf B}$, and denote by $f_{\bf Y}:=\frac{dp_{\bf Y}}{d\mu_W}$ the \textsc{Radon-Nikodym} of the absolutely continuous probability $p_{\bf Y}$ with respect to the \textsc{Wiener} measure. Moreover, assume that $f_{\bf Y}\in L^2(\mu_W)$. Then,  we have
\begin{equation} \rho_{f_{\bf Y}}(J([0,t])) \leq E_{\mathds{P}}\left[ \exp\left( \int_0^t \dot{\widetilde{u}}_s dY_s -\frac{1}{2}\int_0^t \dot{\widetilde{u}}_s^2 ds\right)\right], \label{eqmonf} \end{equation}
with equality, if and only if, there exists a $\mathcal{B}(C([0,1],\mathbb{R}))^{p_{\bf Y}}/\mathcal{B}(C[0,1],\mathbb{R})-$  measurable function $${\bf H}^{u} : C([0,1], \mathbb{R}) \to  C([0,1],\mathbb{R})$$ such that ${\bf H}^{u}_s$ is $\mathcal{F}_s^{p_{\bf Y}}-$ measurable, where ${\bf H}^{u}_s:= W_s\circ  {\bf H}^{u}$, $\forall s\in [0,t]$,  and 
\begin{equation}B_s = {H}^{\bf u}_s({\bf Y}),\ \forall s\in [0,t], \mathbb{P}-a.s..\label{Bequt}\end{equation} Moreover, still in this case,  for any ${\bf h}\in H^1$, and $t\in[0,1]$, $\widehat{h}_s:= H_s^u({\bf Y^h})$, $\mathbb{P}-a.s.$ provides an unbiased estimator for $h_s$ for any $s\in[0,t]$, that is,\begin{equation}h_s = E_{\mathbb{P}}\left[\widehat{h}_s\right], \  \forall s\in[0,t],\label{reconstructf}\end{equation} where  ${\bf Y^h}$ is given by~(\ref{YDEF00}), with ${\bf X} := {\bf B} +{\bf h}$. 

\end{theorem}
\nproof
Lemma~\ref{Lemme1} and Proposition~\ref{Prop1} yield $$\rho_{f_{\bf Y} }(J([0,t])) = E_{\mathbb{P}}\left[ \frac{1}{E_{\mathds{P}}\left[\mathcal{E}(-\delta^{\bf B}\pi_t \widetilde{u})  \middle|   \mathcal{G}_{[0,t]}^Y \right]}\right],$$  so that~(\ref{eqmonf}) follows from the conditional \textsc{Jensen} inequality, by taking into account that $Y_s = B_s +\widetilde{u}_s$, $\forall s\in [0,t]$, outside some $\mathbb{P}-$ negligible set. Since the inverse function is strictly convex on $(0,+\infty)$, and the $\sigma-$ field $\mathcal{G}_{[0,t]}^Y$ is $\mathbb{P}-$ complete, the equality occurs in~(\ref{eqmonf}), if and only if,  $\mathcal{E}(-\delta^{\bf B}\pi_t \widetilde{u})$ is $ \mathcal{G}_{[0,t]}^Y-$ measurable. From the proof of Lemma~\ref{Lemme1} (see~(\ref{cmapong})), the latter condition is also equivalent to the following causal \textsc{Monge-Amp\`{e}re} equation  \begin{equation} E_{p_{\bf B} }\left[\frac{dp_{\bf Y}}{dp_{\bf B} }\middle| \mathcal{F}_t^\nu\right]\circ {\bf Y}   \mathcal{E}(-\delta^{\bf B}\pi_t \widetilde{u}) = 1, \ \mathbb{P}-a.s..\label{CMAE} \end{equation} Then, it is routine to check that the latter is indeed equivalent to the existence of ${\bf H}^{{\textbf{u}}}$, which satisfies the assertions of the claim, by applying same methods as Proposition 3.1. of \cite{FEYEL2}.  However, since the equality is assumed here at time $t$, for the reader's convenience, and for the paper to be enough self-contained, we fix details of the proof for the sake of completeness, and recall why~(\ref{CMAE}) is indeed equivalent to the existence of ${\bf H^u}$ as stated, which in particular satisfies to~(\ref{Bequt}).  First assuming the existence of ${\bf H}^u$ as stated, since ${\bf Y} ={\bf F}^u\circ {\bf B}$ entails $\mathcal{G}_{[0,s]}^{\bf Y} \subset \mathcal{G}_{[0,s]}^{\bf B} $,  we first obtain  $\mathcal{G}_{[0,s]}^{{\bf B} } = \mathcal{G}_{[0,s]}^{\bf Y}$, $\forall s\in [0,t]$. As a consequence, since $(\widetilde{u}_s)$ is $(\mathcal{G}_{[0,s]}^{\bf B} )-$ adapted,  so that together with~(\ref{INTNOMPE2}), \textsc{It\^{o}}'s Lemma ensures that $(\mathcal{E}(-\delta^{\bf B} {\bf \pi_s \widetilde{u}}))_{s\in[0,t]}$ is a $(\mathcal{G}_{[0,s]}^{\bf B} )_{s\in[0,t]}-$ martingale on $(\Omega, \mathcal{A}, \mathbb{P})$, we get \begin{equation} E_{\mathbb{P}}[\mathcal{E}(-\delta^{\bf B} \pi_t {\bf \widetilde{u}})| \mathcal{G}_{[0,t]}^{\bf Y} ] = E_{\mathbb{P}}[\mathcal{E}(-\delta^{\bf B} \pi_t {\bf \widetilde{u})}| \mathcal{G}_{[0,t]}^{\bf B}] = \mathcal{E}(-\delta^{\bf B}  \pi_t {\bf \widetilde{u}}), \ \mathbb{P}-a.s., \label{MMARTEMP} \end{equation}  while similarly as above (see Remark~\ref{remark1M}), the \textsc{Cameron-Martin-Girsanov} theorem yields \begin{equation} E_{p_{\bf B} }\left[\frac{dp_{\bf Y} }{dp_{\bf B} }\middle| \mathcal{F}_t^{p_{\bf B}}\right]\circ {\bf Y}  E_{\mathbb{P}}[\mathcal{E}(-\delta^{\bf B} {\bf \widetilde{u}})| \mathcal{G}_{[0,t]}^{\bf Y} ]     = 1,\ \mathbb{P}-a.s.. \label{MMARTEMP2} \end{equation} Substituting (\ref{MMARTEMP}) into~(\ref{MMARTEMP2}), we obtain~(\ref{CMAE}).

 On the converse, assume now that~(\ref{CMAE}) holds, and recall that since the \textsc{Cameron-Martin-Girsanov} theorem ensures that $p_{\bf Y} $ is absolutely continuous with 
 respect to $p_{\bf B} $,  it is known (see \cite{KAILATH}, see also \cite{DUNCAN}, or \cite{17.}) that there exists a unique \begin{equation}b^{p_{\bf Y} }:=\int_0^. v^{p_{\bf Y}}ds \in 
 L^2_a(p_{\bf Y} , H^1) \label{bpydef} ,\end{equation} which belongs to $L^2_a(p_{\bf Y},H^1)$ since ${\bf \widetilde{u}}\in L^2_a(\mathbb{P}, H^1)$ (see \cite{BZZ} formula $(75)$ 
 and $(76)$ and see also~\cite{17.}), 
 such that setting ${\bf H}^u:= \mathds{I}_{C([0,1], \mathbb{R})}- b^{p_{\bf Y} }$, the stochastic process $(H^{\bf u}_s)_{s\in[0,1]}$ is a $(\mathcal{F}_s^{p_{\bf Y} })-$ brownian 
 motion on $( C([0,1],\mathbb{R}),$ $\mathcal{B}( C([0,1],\mathbb{R})),$ $p_{\bf Y} )$,  where $H^{\bf u}_s:= W_s\circ {\bf H^u}$, $\forall s\in [0,1]$, and
  \begin{equation}E_{p_{\bf B}}\left[\frac{dp_{\bf Y}}{dp_{\bf B}}\middle| \mathcal{F}_{\sigma}^{p_{\bf B}}\right]=\frac{dp_{\bf 
 Y}|_{\mathcal{F}_\sigma^{p_{\bf B}}} }{dp_{\bf B}|_{\mathcal{F}_\sigma^{p_B}} } = \exp\left(\int_0^\sigma v_s^{p_{\bf Y} }  dW_s - \frac{1}{2}\int_0^\sigma (v_s^{p_{\bf Y} } )^2 ds\right), \ p_{\bf Y}-a.s. \label{vsnudefthmeq} ,\end{equation}  $\forall \sigma \in [0,1]$. With this notation, the definitions 
 easily yield $$b^{p_{\bf Y} }\circ {\bf Y}  = \int_0^.E_{\mathbb{P}}\left[\dot{\widetilde{u}}_s| \mathcal{G}_{[0,s]}^{\bf Y} \right]ds, \ \mathbb{P}-a.s.,$$ where the r.h.t denotes the dual 
 predictable projection of $(\widetilde{u}_s)$ on  $(\mathcal{G}_{[0,s]}^{\bf Y} )_{s\in[0,t]}$, as it is well known (for instance see Proposition 3.1. of  \cite{INTRINSIC} for a direct proof).   
 Hence, taking the logarithm and then the expectation under $\mathbb{P}$ in both terms of~(\ref{CMAE}) we get $$(\pi_tb^{p_{\bf Y} })\circ {\bf Y}  =  \int_0^.1_{[0,t]}
 (s)E_{\mathbb{P}}\left[\dot{\widetilde{u}}_s| \mathcal{G}_{[0,s]}^{\bf Y} \right]ds= \pi_t\widetilde{u}, \ \mathbb{P}-a.s.,$$ where we used that since the $L^2(\mathbb{P}, H^1)$ norm of $\pi_t u$ 
 coincides with the norm of  its projection $ \int_0^.1_{[0,t]}(s)E_{\mathbb{P}}\left[\dot{\widetilde{u}}_s| \mathcal{G}_{[0,s]}^{\bf Y} \right]ds$ on the closed linear subspace of the ${\bf 
 k}:=\int_0^. k_{\widetilde{s}}d\widetilde{s} \in L^2_a(\mathbb{P}, H^1)$ such that $(k_{\widetilde{s}})$ is $(\mathcal{G}_{[0,\widetilde{s}]}^{\bf Y})_{\widetilde{s} \in [0,1]}-$ adapted,  
 both elements of $L^2_a(\mathbb{P}, H^1)$ necessarily coincide. Thus, we get \begin{equation} {\bf H}_s^u\circ {\bf F}^u\circ {\bf B} = B_s, \ \forall s\leq t, \ \mathbb{P} -a.s., 
 \label{eqstarmpl}\end{equation} so that, taking into account once again that ${\bf Y} ={\bf F^u}\circ {\bf B} $, $\mathbb{P}-a.s.$, we obtain~(\ref{Bequt}). 

We now turn to the proof of~(\ref{reconstructf}). Thus, assume that we are in the equality case in the above inequality, and still denote by ${\bf H}^u$ the previous map. We then take 
${\bf h}\in H^1$, and denote by ${\bf Y}^{h}:={\bf F}^u({\bf B}+{\bf h})$.  First observe that $${\bf Y}^{h}={\bf F}^u\circ ({\bf B}+{\bf h})={\bf F}^u\circ {\bf T_h}\circ {\bf B},$$ where ${\bf 
T_h}(\omega):={\bf \omega} + {\bf h}\in C([0,1], \mathbb{R})$, $\forall \omega\in C([0,1],\mathbb{R})$. Therefore, taking into account that $p_{T_h({\bf B})}\sim \mu_W$, which follows 
from the \textsc{Cameron-Martin-Segal} theorem (see \cite{MALLIAVIN}, \cite{StroockPTAV}, Theorem 2.2. p. 339 of \cite{YORBOOK}), from~(\ref{eqstarmpl}) we get $$H_s^{\bf u}
\circ {\bf F^u} = W_s, \forall s\leq t, p_{\bf T_h(B)}-a.s.,$$ from which we deduce that  $$H_s^u \circ {\bf Y}^{\bf h} = H_s^{\bf u}\circ {\bf F^u}\circ {\bf T_h}\circ {\bf B} = B_s + h_s, 
\forall s\leq t, \mu_W-a.s..$$ By taking the expectation under $\mathbb{P}$, and since $B_0=0$,~(\ref{reconstructf})  follows. 

\nqed

\section{Second application :  the innovation problem of filtering}
\label{Section5}

The following result may be seen as complementary to the approaches unrolled  in  \cite{ALLINGERM}, \cite{UZ07},  \cite{USTUNEL}, \cite{USTUNEL2} and  \cite{ABSAPI} ; for related topics, see also  \cite{LOCINVSTO} or \cite{RLASU}. It yields a criterion for the innovation conjecture to hold which is based on the \textsc{Tsirelon} spectral measure.

\begin{corollary}\label{Corollary1} 
Let $t\in [0,1]$ and let ${\bf u}\in L^2_a(\mu_W, H^1)$ be such that~(\ref{INTNOMPE2}) holds.  Further define ${\bf Y}$ by~(\ref{YDEF}),  set ${\bf \widetilde{u}}:={\bf Y} -{\bf B}$,  
denote by $f_{\bf Y}:=\frac{dp_{\bf Y}}{d\mu_W}$ the \textsc{Radon-Nikodym} of the absolutely continuous probability $p_{\bf Y}$ with respect to the \textsc{Wiener} measure $\mu_W
$, and assume that $f_{\bf Y}\in L^2(\mu_W)$.  Then,  the innovation conjecture is satisfied by $(Y_s)$ until time $t\in[0,1]$, that is, \begin{equation}(\mathcal{G}_{[0,s]}^{{\bf 
B^{\widetilde{u}}}})_{s\in[0,t]}= (\mathcal{G}_{[0,s]}^{\bf Y})_{s\in[0,t]}, \label{inovconject1} \end{equation} if and only if, there exists a ${\bf \xi}\in L^2_a(\mu_W, H^1)$ which meets the 
following assumptions :
\begin{enumerate}[(i)]
\item  $p_{\bf Y}|_{\mathcal{F}_t^{\mu_W}} =p_{\bf Z}|_{\mathcal{F}_t^{\mu_W}},$ where ${\bf Z}:={\bf F}^{\bf \xi}({\bf B})$, with ${\bf F}^{\bf \xi}:=\mathds{I}_{C([0,1],\mathbb{R})}+ {\bf \xi}$.
\item $E_{\mathbb{P}}\left[ \mathcal{E}(- \delta^{\bf B} \pi_t {\bf \widetilde{\xi}}) \right] =1$, where ${\bf \widetilde{\xi}} :={\bf \xi} \circ {\bf B} = {\bf Z}- {\bf B} \in L^2_a(\mathds{P}, H^1)$.
\item  The following formula holds : \begin{equation}\label{eqxi1}\rho_{f_{\bf Y}}(J([0,t])) = E_{\mathds{P}}\left[ \exp\left( \int_0^t \dot{\widetilde{\xi}}_s dZ_s -\frac{1}{2}\int_0^t 
\dot{\widetilde{\xi}}_s^2 ds\right)\right],\end{equation} where $$J([0,t]):= \{K\in Comp([0,1]) : K\subset [0,t]
\}.$$  
\end{enumerate}
\end{corollary}
\nproof
First assuming the existence of a ${\bf \xi}\in L^2_a(\mu_W,H^1)$ as above, notice that  $(i)$ implies $$E_{\mu_W}\left[\frac{d p_{\bf Y}}{d\mu_W}\middle|\mathcal{F}_{[0,t]}^{\mu_W}
\right] =E_{\mu_W}\left[\frac{d p_{\bf Z}}{d\mu_W} \middle| \mathcal{F}_{[0,t]}^{\mu_W}\right],\ \mu_W-a.s.,$$ so that $$\rho_{f_{\bf Y}}(J([0,t]))=\rho_{f_{\bf Z}}(J([0,t])),$$ where 
$f_{\bf Z}:= \frac{d p_{\bf Z}}{d
\mu_W}$. Thus, by $(ii)$ and $(iii)$, Theorem~\ref{Theorem1TSMF} applied to ${\bf Z}$ yields  the existence of ${\bf H}^{\bf \xi}$ such that ${\bf H}^\xi_s\circ {\bf Z} = B_s$, $\forall s
\leq t$, outside some $\mathbb{P}-$ null set. Since $Z_s= {\bf F}^\xi_s\circ {\bf B}$, for all $s\leq t$, $\mathbb{P}-a.s.$, it entails \begin{equation} H^{\bf \xi}_s\circ {\bf F}^\xi = W_s, 
\forall s\leq t,  \ \mu_W-a.s. \label{corollarystarad} ,\end{equation} which notably implies that $(H^{\bf \xi}_s)_{s\in[0,t]}$ is a canonical Brownian motion until time $t$, on the 
canonical probability space endowed with the probability $p_Z$.

Let $A := \{\omega \in C([0,1],\mathbb{R}) : H^{\bf \xi}_s\circ {\bf F}^\xi = W_s, \forall s\leq t \},$ $\widetilde{A}:= \{\omega \in C([0,1],\mathbb{R}) : F^{\bf \xi}_s \circ {\bf H}^{\bf \xi}\circ 
{\bf F}^{\bf \xi} = F^{\bf \xi}_s, \forall s\leq t \},$ and $C=\{\omega \in C([0,1],\mathbb{R}) :  F^{\bf \xi}_s\circ {\bf H}^{\bf \xi} = W_s, \forall s\leq t \}$.
 Since $A\subset \widetilde{A}$ and $\mu_W(A)=1$, by noticing that $C\in \mathcal{F}_t^{p_{\bf Z}}$, we get $$1\geq p_{\bf Z}(C) ={\bf F^\xi}_\star\mu_W(C)= \mu_W(\widetilde{A}) 
 \geq \mu_W(A)=1,$$ ${\bf F^\xi}_\star{\mu_W}$ denoting the pushforward of the probability $\mu_W$ by ${\bf F^\xi}$, so that \begin{equation}F_s^{{\bf \xi}}\circ {\bf H}^{\bf \xi} 
 =W_s, \ \forall s\in[0,t],  \ p_{\bf Z}-a.s.. \label{fxiinv2}  \end{equation} From this, we deduce that $H_s^{\bf \xi} = W_s-\xi_s\circ {\bf H^\xi}$, $\forall s\leq t$, $p_Z|_{\mathcal{F}_t}-a.s.
 $.  Since $p_{\bf Y}|_{\mathcal{F}_t^{\mu_W}}=  p_{\bf Z}|_{\mathcal{F}_t^{\mu_W}}$ and $p_{\bf Y}\sim \mu_W$ (equivalent), 
 the uniqueness of the semi-martingale decomposition of $(W_s)_{s\in[0,t]}$ entails,  $\xi_s \circ {\bf H^\xi} = b_s^{p_{\bf Y}}$,  $\forall s\leq t$ $p_{\bf Y}-a.s.$ (and then $p_{\bf B}-
 a.s.$), where ${\bf b^{p_{\bf Y}}}:=\int_0^. v_s^{p_{\bf Y}}ds$ is as in~(\ref{bpydef}), so that in particular  $(v_s^{p_{\bf Y}})$ satisfies to~(\ref{vsnudefthmeq}) (finite variation part), and 
 $B^{\widetilde{u}}_s =H_s^{\bf \xi}\circ {\bf Y}$, $\forall s\leq t$, $\mathbb{P}-a.s.$, which easily follows from the identification of the martingale parts together with the definition of the 
 innovation process $(B^{\widetilde{u}}_s)$. Hence~(\ref{fxiinv2}) yields $F^{\xi}_s\circ {\bf B^{\widetilde{u}}} = Y_s, $ $\forall s\in [0,t]$, $\mathbb{P}-a.s.$, and therefore $\mathcal{G}_{[0,s]}^{\bf Y}\subset \mathcal{G}_{[0,s]}^{{\bf B^{\widetilde{u}}}}$, $\forall s\in [0,t]$. Since $(B^{\widetilde{u}}_s)$ is known to be a $(\mathcal{G}_{[0,s]}^{\bf Y})_{s\in[0,1]}-$ 
 brownian motion,  this proves~(\ref{inovconject1}).

 Conversely, we now assume that~(\ref{inovconject1}) holds, and define the continuous stochastic processes $(\widetilde{Y}_s)$ and $(\widetilde{X}_s)$ by $$\widetilde{Y}_s = 1_{[0,t)}(s) Y_s + 1_{[t,1]}(s)(Y_t+ B_s-B_t), \ \forall s\in[0,1], \mathbb{P}-a.s.,$$ and $$\widetilde{X}_s= 1_{[0,t)}(s)B_s^{\bf \widetilde{u}} + 1_{[t,1]}(s) (B_t^{ \bf \widetilde{u}} +B_s-B_t) \ \forall s\in[0,1], \mathbb{P}-a.s.$$
 Then,~(\ref{inovconject1}) implies $$(\mathcal{G}_{[0,s]}^{{\bf \widetilde{Y}}})_{s\in[0,1]}= (\mathcal{G}_{[0,s]}^{\widetilde{X}})_{s\in[0,1]}.$$ 
This ensures  the existence of  two maps ${\bf \widetilde{I}},{\bf \widetilde{J}} : C([0,1],\mathbb{R}) \to C([0,1],\mathbb{R})$ such that ${\bf \widetilde{I}}$ is $
 \mathcal{B}^{p_{\widetilde{X}}}(C([0,1],\mathbb{R}))/ \mathcal{B}(C([0,1],\mathbb{R}))$ - measurable, while ${\bf \widetilde{J}}$ is $\mathcal{B}
 ^{p_{{\bf \widetilde{Y}}}}(C([0,1],\mathbb{R}))/$ $\mathcal{B}(C([0,1],\mathbb{R}))$ - measurable, ${\bf \widetilde{J}}\circ {\bf \widetilde{I}} = \mathds{I}_{C([0,1],\mathbb{R})}$, 
 $p_{\widetilde{X}}-a.s.$,  ${\bf \widetilde{I}}\circ {\bf \widetilde{J}} = \mathds{I}_{C([0,1],\mathbb{R})}$, $p_{{\bf \widetilde{Y}}}-a.s.$, ${\bf \widetilde{I}}$ (resp. ${\bf \widetilde{J}}$) is $(\mathcal{F}_t^{p_{\widetilde{X}}})-$ (resp.  $(\mathcal{F}_t^{p_{{\bf \widetilde{Y}}}})-$) adapted, and ${\bf \widetilde{Y}}= {\bf \widetilde{I}}\circ {\bf \widetilde{X}}$ (resp. $
 {\bf \widetilde{X}}= {\bf \widetilde{J}}\circ {\bf \widetilde{Y}}$), $\mathbb{P}-a.s.$ (for instance, see Proposition~1.3. of \cite{INTRINSIC}).   Hence, we obtain \begin{equation}Y_s=\widetilde{Y}_s = \widetilde{I}_s\circ {\bf \widetilde{X}}
 =\widetilde{I}_s\circ {\bf B^{\widetilde{u}}}, \ \forall s\leq t, \ \mathbb{P}-a.s., \label{Ydefmplsep} \end{equation} and $B_s^{\widetilde{u}}=\widetilde{X}_s = \widetilde{J}_s\circ 
 {\bf \widetilde{Y}}=\widetilde{J}_s\circ {\bf Y}$, $\forall s\leq t$, $\mathbb{P}-a.s.$. On the other hand, since $B^{\widetilde{u}}= (\mathds{I}_{C([0,1],\mathbb{R})}-b^{p_{\bf Y}})\circ {\bf Y}$, 
 where  $b^{p_{\bf Y}}$ is defined by~(\ref{bpydef}), we get $\widetilde{J}_s\circ {\bf Y} = Y_s -b_s^{p_{\bf Y}}\circ {\bf Y}$, $\forall s\in [0,t]$, $\mathbb{P}-a.s.$, so that  $\widetilde{J}_s =W_s - b_s^{p_{\bf Y}}$, $\forall s\in [0,t]$, $p_{\bf Y}|_{\mathcal{F}_t^{p_{\bf Y}}}-$ a.s., which implies that $W_s=\widetilde{J}_s\circ {\bf \widetilde{I}} = \widetilde{I}_s-\xi_s$, $\forall s\leq t$, $\mu_W-$a.s., where we have defined $\xi := (\pi_t b^{p_{\bf Y}})\circ {\bf \widetilde{I}} = \int_0^. 1_{[0,t]}(s)v_s^{p_{\bf Y}}\circ {\bf \widetilde{I}} ds$. Further notice that $p_{\bf \widetilde{X}}=\mu_W$, so that ${\bf \xi}\in L^2_a(\mu_W,H^1)$ follows from $b^{p_{\bf Y}}\in L^2_a(p_{\bf Y},H^1)$. Then, we define ${\bf F^{\xi}}:= \mathds{I}_{C([0,1],\mathbb{R})} + {\bf \xi}$, ${\bf \widetilde{\xi}}:= {\bf \xi}\circ {\bf B}$, and ${\bf Z}:={\bf F^{\xi}}({\bf B}) ={\bf B} +{\bf \widetilde{\xi}}$. By definition, we obtain $Z_s =\widetilde{I}_s({\bf B})$, $\forall s\leq t$, $\mathbb{P}-a.s.$, so that in particular, $\widetilde{\xi}_s= b_s^{p_{\bf Y}}\circ {\bf Z}$, $\forall s\leq t$, $\mathbb{P}-a.s.$, while we obtain $\widetilde{I}_s =W_s +\xi_s$, $\forall s\leq t$, $\mu_W-a.s.$. Thus, taking into account that $p_{\bf B^{{\widetilde{u}}}}= p_{\bf B}$,~(\ref{Ydefmplsep}) entails that $p_{\bf Z}|_{\mathcal{F}_t^{\mu_W}}= p_{{\bf Y}}|_{\mathcal{F}_t^{\mu_W}}$, while $(\mathcal{G}_{[0,s]}^{\bf Z})_{s\in[0,t]} \subset (\mathcal{G}_{[0,s]}^{\bf B})_{s\in[0,t]}$. On the other hand, due to the \textsc{Camron-Martin-Girsanov} theorem, ~(\ref{INTNOMPE2})  
  entails $p_{\bf Y}|_{\mathcal{F}_t^{\mu_W}}\sim \mu_W |_{\mathcal{F}_t^{\mu_W}}$, so that  (\ref{vsnudefthmeq}) yields  \begin{equation}\frac{{dp_{\bf Y}}|_{\mathcal{F}_t^{\mu_W}} }{d{p_B}|_{\mathcal{F}_t^{\mu_W}}} = \exp\left(\int_0^t v_s^{p_{\bf Y}} dW_s - \frac{1}{2} \int_0^t (v_s^{p_{\bf Y}})^2 ds  \right), \ {\mu_W}|_{\mathcal{F}_t^{\mu_W}}-a.s. \label{eqinme}  . \end{equation}   Therefore, we get 
 \begin{eqnarray*}
 E_{\mathbb{P}}\left[\mathcal{E}(-\delta^{B}\pi_t \widetilde{\xi})\right] & = & E_{\mathbb{P}}\left[ \exp\left(-\int_0^t v_s^{p_{\bf Y}}\circ {\bf Z} dB_s -\frac{1}{2} \int_0^t (v_s^{p_{\bf 
 Y}}\circ {\bf Z})^2 ds    \right)\right] \\  & = & E_{\mu_W}\left[ \exp\left(-\int_0^t v_s^{p_{\bf Y}}\circ \widetilde{I}_s d\widetilde{I}_s + \frac{1}{2} \int_0^t (v_s^{p_{\bf Y}}\circ 
 \widetilde{I}_s)^2 ds  \right)\right]  \\ & = &  E_{p_{\bf Z }}\left[\exp\left(-\int_0^t v_s^{p_{\bf Y}} dW_s + \frac{1}{2} \int_0^t (v_s^{p_{\bf Y}})^2 ds  \right)\right]    \\ & = &  E_{p_{\bf Y }}
 \left[\exp\left(-\int_0^t v_s^{p_{\bf Y}} dW_s + \frac{1}{2} \int_0^t (v_s^{p_{\bf Y}})^2 ds  \right)\right]    \\ & = &  E_{\mu_W|_{\mathcal{F}_t^{\mu_W}}}\left[ \frac{d p_{\bf Y}|_{\mathcal{F}_t^{\mu_W}}}{d \mu_W|_{\mathcal{F}_t^{\mu_W}}} \exp\left(-\int_0^t v_s^{p_{\bf 
 Y}} dW_s + \frac{1}{2} \int_0^t (v_s^{p_{\bf Y}})^2 ds  \right)\right]   \\ & = & 1, \end{eqnarray*} where  the last inequality follows from~(\ref{eqinme}). Thus,  ${\bf \xi}$ satisfies to $(i)$ and $(ii)$. Finally, notice that for $s\leq t$, we have $
 \widetilde{J}_s\circ {\bf Z} = \widetilde{J}_s\circ {\bf \widetilde{I}}\circ {\bf B} = W_s\circ {\bf B}= B_s$, $\forall s\leq t$, $\mathbb{P}-a.s.$, so that  $(\mathcal{G}_{[0,s]}^{\bf B})_{s
 \in[0,t]} = (\mathcal{G}_{[0,s]}^{\bf Z})_{s\in[0,t]}$. As a consequence, we are in the equality case of Theorem~\ref{Theorem1TSMF} applied to ${\bf \xi}$, from which we obtain that 
 the latter satisfies to $(iii)$. This achieves the proof.
\nqed

\section{Probabilistic perspective by normalization, spectral random sets,  and further developments}
\label{Section6}
To clarify the statements in view of applications to the above problems of filtering, we did not require the \textsc{Tsirelson} spectral measures  to be necessarily normalized. However, taking a finite measure of finite mass one, as it is usual in stochastic analysis, notably increases the insight on this kind of studies, as it is the case in Example~\ref{Example1} and Example~\ref{Example2} above. Moreover, we notice that the spectral measures investigated above were determined by the laws of the associated processes. This motivates to introduce the following notation to shorten the statement of Corollary~\ref{emptynot} below :

\begin{definition} \label{ropfmr}
Let $\nu\in {\bf M}_1(C([0,1],\mathbb{R}))$ be a \textsc{Borel} probability measure which is absolutely continuous with respect to the \textsc{Wiener} measure ($\nu<< \mu_W$), whose \textsc{Radon-Nikodym} derivative $f_{\bf Y}:= \frac{d\nu}{d\mu_W}$ is square-integrable (i.e. $E_{\mu_W}\left[(\frac{d\nu}{d\mu_W})^2\right] <+\infty$). Then, we denote by  $ \mathbb{P}_{\nu}$ the normalized \textsc{Tsirelson} spectral measure associated with the probability $\nu$,  which is the \textsc{Borel} probability measure $\mathbb{P}_\nu$ on $Comp([0,1])$ 
given by \begin{equation}\mathbb{P}_\nu(\mathcal{K}) := \frac{\rho_{\frac{d\nu}{d\mu_W}}(\mathcal{K})}{\rho_{{\frac{d\nu}{d\mu_W}}}(\text{Comp}([0,1]))},\ \forall \mathcal{K}\in \mathcal{B}(\text{Comp}([0,1])), \label{spectprobdef}\end{equation} i.e.,  $$\mathbb{P}_\nu := \rho_{g_\nu} \in {\bf M}_1(Comp([0,1])),$$ where we have set $$g_\nu := {\frac{\frac{d\nu}{d\mu_W}}{\sqrt{E_{\mu_W}[({\frac{d\nu}{d\mu_W})^2]}}}},$$ and where $\rho_{g_\nu}$ denotes the \textsc{Tsirelson} spectral measure (see Definition~(\ref{Tsimesdef})) of $\mu_W$ with respect to $g_\nu\in L^2(\mu_W)$.
\end{definition}

Notice that the latter is well defined since  from  \textsc{Jensen}'s inequality, it follows that  $\|\frac{d\nu}{d\mu_W}\|_{L^2(\mu_W)}>0$, for any $\nu<<\mu_W$ (absolutely continuous). In this way, under the assumptions of Definition~\ref{ropfmr}, any such  $\nu$ is associated to a random compact set ; in the limit case where $\nu=\mu_W$ this compact set is a.s. empty.  Then, we obtain the following ersatz for part of Theorem~\ref{Theorem1TSMF}, which notably provides a probabilistic insight on the previous criterion of existence of a realizable inverse to a real-time noise filter, thanks to some spectral random set  :

\begin{corollary}\label{emptynot}
Given  ${\bf u}\in L^2_a(\mu_W, H^1)$ which satisfies to~(\ref{INTNOMPE2}) with $t=1$, further define ${\bf Y}$ by~(\ref{YDEF}),  and set ${\bf \widetilde{u}}:={\bf Y} -{\bf B}$. Further assume that we have $\frac{d\nu}{d\mu_W}\in L^2(\mu_W)$, the latter denoting the \textsc{Radon-Nikodym} derivative of the  absolutely continuous probability $\nu:=p_{\bf Y}$ with respect to the  \textsc{Wiener} measure $\mu_W$.Then,  the following lower bound holds  $$\mathbb{P}_{\nu}(\{\emptyset\}) \geq \frac{1}{E_{\mathds{P}}\left[ \exp\left( \int_0^1 \dot{\widetilde{u}}_s dY_s -\frac{1}{2}\int_0^1 \dot{\widetilde{u}}_s^2 ds\right)\right]},$$  $\mathbb{P}_{\nu}$ denoting the normalized \textsc{Tsirelson} spectral measure~(\ref{spectprobdef}) of the probability $p_{\bf Y}$.  Moreover, there exists ${\bf H}^{\bf u} : C([0,1], \mathbb{R}) \to  C([0,1],\mathbb{R})$ which satisfies to the hypothesis of Theorem~\ref{Theorem1TSMF} with $t=1$, if and only if, $$\mathbb{P}_{\nu}(\{\emptyset\}) = \frac{1}{E_{\mathds{P}}\left[ \exp\left( \int_0^1 \dot{\widetilde{u}}_s dY_s -\frac{1}{2}\int_0^1 \dot{\widetilde{u}}_s^2 ds\right)\right]}.$$ \end{corollary}
\nproof
Since $Comp([0,1])= J([0,1])$ (see~(\ref{JDEFPME}) above), from Definition~\ref{ropfmr}, we get $\mathbb{P}_{\nu}(\{\emptyset\}) = \frac{\rho_{f_{\bf Y}}(\{\emptyset\})}{\rho_{f_{\bf Y}}(J([0,1]))},$ where $f_{\bf Y}:=\frac{dp_{\bf Y}}{d\mu_W}$. On the other hand,~(\ref{tsimot}) implies $$\rho_{f_{\bf Y}}(\{\emptyset\}) = \left(E_{\mu_W}\left[f_{\bf Y}\right]\right)^2  = 1,$$ since $f_{\bf Y}$ is the \textsc{Radon-Nikodym} derivative of a probability measure with respect to the probability $\mu_W$.  Thus, the result follows from Theorem~\ref{Theorem1TSMF}.
\nqed

\begin{remark}
Within this viewpoint, the above results can be extended in several fashions to a class of \textsc{Borel} probability measures absolutely continuous with respect to the \textsc{Wiener} measure whose \textsc{Radon-Nikodym} derivatives are also square integrable, which will be done in some forthcoming works.
\end{remark}

.


\begin{thebibliography}{}


\bibitem{AirMal} \textsc{Airault, H.}, \textsc{Malliavin, P.} : \textit{Int\'{e}gration et analyse de Fourier, Probabilit\'{e}s et analyse gaussienne.}  $2^e$ \'{e}dition revue et augment\'{e}e. MASSON, Paris Milan Barcelone (1994)


\bibitem{ALLINGERM} \textsc{Allinger, D.}, \textsc{Mitter S.K.} :  New results on the innovations problem for nonlinear filtering. \textit{Stochastics} 4, no.4, 339$-$348 (1980)

\bibitem{BAIN} \textsc{Bain, A.}, \textsc{Crisan, D.}: \textit{Fundamentals of Stochastic Filtering.} Stochastic Modeling And Applied Probability 60. Springer (2009)

 \bibitem{BZZ} \textsc{Binia, J.}, \textsc{Zakai, M.}, \textsc{Ziv, J.} : On the epsilon-entropy and the rate-distortion function of certain non-Gaussian processes. IEEE {Transactions on Information Theory}, vol. IT$-$20, no. 4, pp. 517-524 (1974).

\bibitem{34.} \textsc{Cameron, R.H.}, \textsc{Martin, W.T.} : Transformation of Wiener integral under translation. \textit{Ann. Math.} 45 (1944)
\bibitem{35.} \textsc{Cameron, R.H.}, \textsc{Martin, W.T.}: The transformation of Wiener integrals by nonlinear transformation. \textit{Trans. Am. Math. Soc.} 66  253-283 (1949)



\bibitem{CONT} \textsc{Cont, R.}, \textsc{Tankov, P.} : \textit{Financial Modelling with Jump Processes} Chapman and Hall/CRC Boca Ra ton London New York Washington, D.C. (2004)




\bibitem{DM1}  \textsc{Dellacherie, C.} and \textsc{Meyer, P. A.} : \textit{Probabilit\'es et Potentiel  Ch. 1 \`a 4}. Paris, Hermann. (1975)

\bibitem{DM5}  \textsc{Dellacherie, C.}, \textsc{Maisonneuve, B.}, \textsc{Meyer, P. A.} :  \textit{Probabilit\'es et Potentiel} V, processus de Markov,  Ch. 17 \`a 24. Paris, Hermann. (1992)


\bibitem{DUNCAN}\textsc{Duncan, T.E.} : Evaluation of likelihood functions. \textit{Inform. Contr.}, vol. 13, pp.62-74, (1968)


\bibitem{FEYEL1} \textsc{Feyel, D.}, \textsc{\"{U}st\"{u}nel, A.S.} : Monge-Kantorovitch Measure Transportation and Monge-Amp\`{e}re Equation on Wiener Space. \textit{Probab. Theory Relat. Fields} 128, 347?385 (2004).

\bibitem{FEYEL2} \textsc{Feyel, D.},  \textsc{\"{U}st\"{u}nel, A.S.},  \textsc{Zakai, M.} : The realization of positive random variables via absolutely continuous transformations of measure on Wiener space.  \textit{Probab. Surveys} 3 170 - 205 (2006).


\bibitem{17.} \textsc{F\"{o}llmer, H.} :  Random fields and diffusion processes. In: Hennequin PL. (eds) \textit{ \'{E}cole d'\'{e}t\'{e} de Saint Flour XV$-$XVII, 1985$-$1987}. Lecture Notes in Mathematics, vol 1362. Springer, Berlin, Heidelberg  (1988) 

\bibitem{GARBAN} \textsc{Garban, C.} : Oded Schramm's contributions to noise sensitivity. \textit{Ann. Probab.}, Vol. 39, No 5, 1702-1767 (2011)


\bibitem{GIRSANOV} \textsc{Girsanov, I.V.} : On transforming a certain class of stochastic processes by absolutely continuous substitution of measures.  {\em Theory Probab. Appl} {\bf 5} 285-301 (1960)

\bibitem{IKEDA} \textsc{Ikeda, N.}, \textsc{Watanabe, S.} : \textit{Stochastic Differential Equations and Diffusion Processes.} Second Edition. North-Holland Mathematical Library, Volume 24. North Holland Publishing Company, Amsterdam (Kodansha Ltd., Tokyo) (1989)

\bibitem{44.} \textsc{It\^o, K.} : Stochastic integral. \textit{Proc. Imp. Acad. Tokyo}  20, no.~8, 519-524 (1944)



\bibitem{26.} \textsc{Kailath, T.} : Some Extensions of the Innovations Theorems. \textit{B.S.T.J.}, 50, p.1487-1494 (1971)

\bibitem{KZ} \textsc{Kailath, T.}, \textsc{Zakai, M.} : Absolute continuity and Radon-Nikodym derivatives for certain measures relative to Wiener measure. \textit{The Annals of Mathematical Statistics}, Vol. 42, No 1, 130-140 (1971)
\bibitem{KAILATH} \textsc{Kailath, T.} : The structure of Radon-Nikodym derivatives with respect to Wiener and related measures. \textit{The Annals of Mathematical Statistics}, Vol. 42, No 3, 1054-1067 (1971)

\bibitem{KOOPMANS} \textsc{Koopmans, L.H.} : \textit{The spectral analysis of time series}. Academic press.  (1974)


\bibitem{KUO} \textsc{Kuo, H.} : \textit{Gaussian Measures in Banach Spaces} . Lecture Notes in Math.  Springer$-$Verlag Berlin Heidelberg New York  (1975)



\bibitem{INTRINSIC} \textsc{Lassalle, R.}, \textsc{Cruzeiro, A.B.} : An intrinsic calculus of variations for functionals of laws of semi-martingales.  \textit{Stochastic Processes and their Applications}, Volume 129, Issue 10, 3585-3618 (2019) 


\bibitem{LOCINVSTO}  \textsc{Lassalle, R.}  : Local invertibility of adapted shifts on Wiener space, under finite energy condition, \textit{Stochastics} Volume {\bf 85}, 2013 - Issue 6, 987-996.  (2013)

\bibitem{RLASU} \textsc{Lassalle, R.}, \textsc{\"{U}st\"{u}nel, A.S.} :  Local Invertibility of Adapted Shifts on Wiener Space, and Related Topics. In : Viens F., Feng J., Hu Y., Nualart E. (eds)  \textit{Malliavin Calculus and Stochastic Analysis.}   Springer Proceedings in Mathematics \& Statistics, vol. 34. Springer, Boston, MA (2013)





\bibitem{ABSAPI} \textsc{Lassalle, R.} : Invertibility of adapted perturbations of the identity on abstract Wiener space.  In  \textit{Journal of Functional Analysis}, volume  262, Issue 6, 2734$-$2776  (2012)

\bibitem{LEJANRAIMOND} \textsc{Le Jan, Y.},  \textsc{Raimond, O.} : Flows, coalescence and noise. \textit{Ann. Probab.} 32 (2) 1247 - 1315 (2004)

\bibitem{MALLIAVIN} \textsc{Malliavin, P.} :  \textit{Stochastic Analysis.} Springer$-$Verlag Berlin Heidelberg, Grundlehren der mathematischen Wissenschaften, volume 313 (1997)


\bibitem{MEYERF} \textsc{Meyer, P.-A.} : Sur un probl\`{e}me de filtration. \textit{S\'{e}minaire de probabilit\'{e}s de Strasbourg}, Tome 7, pp. 223-247 (1973)

\bibitem{MEYERQUANTUM} \textsc{Meyer, P.-A.} : \textit{Quantum probability for probabilists.} Second edition. Lecture notes in mathematics 1538. Springer-Verlag Berlin Heidelberg (1995)

\bibitem{NOURDIN} \textsc{Nourdin, I.} \and \textsc{Peccati, G.}  \textit{Normal approximations with Malliavin Calculus.} Cambridge University Press.  192. (2012)

\bibitem{PRIVAULT} \textsc{Privault, N.} : \textit{Stochastic analysis in discrete and continuous settings} Lecture notes in mathematics 1982. Springer-Verlag Berlin Heidelberg (2009)

\bibitem{PROTTER} \textsc{Protter, P.E.} : \textit{Stochastic integration and differential equations.} Stochastic Modeling and Applied Probability 21. Second edition, version 2.1. Springer Berlin Heidelberg New York. Second edition, 3rd Printing (2005)

\bibitem{YORBOOK} \textsc{Revuz, D.}, \textsc{Yor, M.}: \textit{Continuous martingales and Brownian motion.} Springer Grundlehren der mathematischen Wissenshaften 293. 3rd edition. Springer Berlin Heidelberg New York (1999) 

\bibitem{StroockPTAV} \textsc{Stroock, D.W.} : \textit{Probability theory, an analytic view.} Cambridge university press, New York. second edition (2011)

\bibitem{Tsirelinvov}\textsc{Tsirelson, B.} : Within and beyond the reach of Brownian innovation, \textit{Documenta Mathematica Extra Volume ICM}, III, 311-326. (1998)
\bibitem{TsirelsonRev}\textsc{Tsirelson, B.} : Nonclassical stochastic flows and continuous products, \textit{Probability Surveys}, 1 (2004), 173-298.

\bibitem{TsirelsonSaintFlour}\textsc{Tsirelson, B.} : Scaling Limit, Noise, Stability. In: Picard, J. (eds) \textit{Lectures on Probability Theory and Statistics. Lecture Notes in Mathematics, vol 1840.} Springer, Berlin, Heidelberg. (2004)


\bibitem{UZ07} \textsc{\"{U}st\"{u}nel, A.S.}, \textsc{Zakai, M.} : Sufficient conditions for the invertibility of adapted perturbations of identity on the Wiener space. \textit{Probab. Theory Relat. Fields} 139, 207-234 (2007).
 \bibitem{USTUNELBOOK}  \textsc{\"{U}st\"{u}nel, A.S.}, \textsc{Zakai, M.} : \textit{Transformation of Measure on Wiener Space} Springer monographs in Mathematics. Springer-Verlag Berlin Heidelberg (2000)


 \bibitem{USTUNEL} \textsc{\"Ust\"unel, A.S.} :  Entropy, invertibility and variational calculus of adapted shifts on Wiener space. \textit{J. Funct. Anal.} 257,  no. 11, 3655-3689 (2009)
 \bibitem{USTUNEL2}\textsc{\"Ust\"unel, A.S.} : Entropic solution of the innovation conjecture of T. Kailath. \textit{Kyoto J. Math.} 55(3): 555-566. (2015)
 
 \bibitem{WARRENWAT}\textsc{Warren, J.}, \textsc{Watanabe, S.} : On Spectra of Noises Associated with Harris Flows.
\textit{Adv. Stud. Pure Math.}, 2004: 351-373 (2004)
 
 \bibitem{Watanabe15}\textsc{Watanabe, S.} : It\^{o}'s theorems on chaos expansions and martingale representations. \textit{Journal of the Mathematical Society of Japan.} Vol. 67, No. 4 (2015)

 \bibitem{WIENER}\textsc{Wiener, N.} : Differential spaces. 131-174 \textit{J. Math. Phys.} 2 (1923)




\end{thebibliography}
\end{document}